\documentclass[12pt]{article}
\usepackage{color,bbold}
\usepackage{amssymb,amsmath,amsthm,graphicx, color, amsfonts,tabularx}
\usepackage[english]{babel}
\usepackage{afterpage}
\usepackage{authblk, relsize}
\usepackage{float}
\usepackage{geometry}
\def\supp{\mathop{\rm supp}\nolimits}

\newtheorem{theorem}{Theorem}[section]

\newtheorem{lemma}[theorem]{Lemma}
\newtheorem{proposition}[theorem]{Proposition}

\newtheorem{definition}[theorem]{Definition}
\newtheorem{remark}[theorem]{Remark}
\newtheorem{example}[theorem]{Example}

\renewenvironment{proof}[1][.]{%
\bigskip\noindent{\bf Proof#1 }}{%
\hfill$\blacksquare$\bigskip}

\begin{document}
\pagestyle{myheadings}
\title{On the Higuchi fractal dimension of invariant measures for countable idempotent iterated function systems}
\author[1]{Elismar R. Oliveira
\thanks{E-mail: elismar.oliveira@ufrgs.br}}
\affil[1]{Universidade Federal do Rio Grande do Sul}

\date{\today}
\maketitle

\begin{abstract}
We investigate the set of invariant idempotent probabilities for countable idempotent iterated function systems (IFS) defined in compact metric spaces. We demonstrate that, with constant weights, there exists a unique invariant idempotent probability. Utilizing Secelean's approach to countable IFSs, we introduce partially finite idempotent IFSs and prove that the sequence of invariant idempotent measures for these systems converges to the invariant measure of the original countable IFS. We then apply these results to approximate such measures with discrete systems, producing, in the one-dimensional case, data series whose Higuchi fractal dimension can be calculated. Finally, we provide numerical approximations for two-dimensional cases and discuss the application of generalized Higuchi dimensions in these scenarios.
\end{abstract}
\vspace {.8cm}



\section{Introduction}\label{sec:intro}
The theory of idempotent probabilities was introduced by Maslov in \cite{KM89}, \cite{LitMas96} and \cite{Kol97} (followed by many other papers and applications in the subsequent years) and further developed by \cite{Rep10} for idempotent measures in compact spaces, \cite{Zar}, \cite{ZAI} for the study of the space of idempotent measures and \cite{BZ2020} for the theory of max-min idempotent measures.

In the present paper our focus is the description of the invariant idempotent probabilities for a  transfer operator associated to a countable iterated function system (CIFS). See \cite{Sec01}, \cite{secelean2002sufficient}, \cite{Secelean2011}, \cite{Secelean2014} and \cite{Secelean2014Invar} for details on CIFSs. The existence of a unique invariant idempotent measure for an IFS was given in \cite{MZ} from a topological point of view and lately in \cite{dCOS20fuzzy} and \cite{ExistInvIdempotent} from metrical point of view. The full characterization of all idempotent invariant measures for a max-plus IFS, in compact metric spaces, is provided by \cite{MO24}. Another key reference is \cite{mazurenko2023invariant} (see also \cite{Sukhorukova_2023}). There we found a generalization of the notion of idempotent measure, called $\ast$-idempotent measure, defined for every triangular norm. They prove existence and uniqueness of invariant $\ast$-idempotent measures for iterated function systems on compact metric spaces.

The paper is organized as follows. In Section~\ref{sec:Preliminaries}, we provide background information on idempotent mathematics and Iterated Function Systems (IFSs). Section~\ref{sec:approx by contract} establishes the existence of the attractor for a countable max-plus IFS through contractivity and defines pointwise approximation by partial attractors. In Section~\ref{sec: fuzzy approach}, we review the fuzzy approach for idempotent invariant measures from \cite{dCOS20fuzzy} to prepare for computing the Higuchi fractal dimension. Finally, Section~\ref{Hichi approximations} focuses on computing the Higuchi fractal dimension using numerical approximations of the partial attractors.

\section{Preliminaries}\label{sec:Preliminaries}
Let us consider the  max-plus semiring $\mathbb{R}_{\max}:=\mathbb{R}\cup \{-\infty\}$ endowed with the operations
\begin{enumerate}
	\item $\oplus: \mathbb{R}_{\max} \times \mathbb{R}_{\max} \to \mathbb{R}_{\max}$, where $a \oplus b :=\max(a,b)$ assuming $a \oplus -\infty:=a$. The max-plus \emph{additive} neutral element is $\mathbb{0} :=-\infty$.
	\item $\odot: \mathbb{R}_{\max} \times \mathbb{R}_{\max} \to \mathbb{R}_{\max}$, where $a \odot b :=a+b$ assuming $a \odot -\infty:=-\infty$. The max-plus \emph{multiplicative} neutral element is $\mathbb{1} :=0$.
\end{enumerate}

We consider also a compact metric space $(X,d)$ and $C(X,\mathbb{R})$, the set of continuous functions from $X$ to $\mathbb{R}$.
We notice that $\mathcal{V}:=(C(X,\mathbb{R}), \oplus, \odot)$ has a natural $\mathbb{R}$-semimoduli (a vectorial space over a semiring) structure:
\begin{enumerate}
	\item $(a \odot f)(x):= a \odot f(x)$, for $a \in \mathbb{R}$ and $f\in C(X,\mathbb{R})$;
	\item $(f \oplus g)(x):=f (x)\oplus g(x)$ for $f, g \in C(X,\mathbb{R})$.
\end{enumerate}

A function $m: C(X,\mathbb{R}) \to \mathbb{R}$ is a max-plus linear  functional if
	\begin{enumerate}
		\item  $m(a \odot f)=a \odot m( f)$, $\forall a \in \mathbb{R}$ and $\forall f \in C(X,\mathbb{R})$(max-plus homogeneity);
		\item  $m(f \oplus g)=m( f) \oplus m( g)$, $\forall f, g \in C(X,\mathbb{R})$ (max-plus additive).
	\end{enumerate}

	We denote by $C_*^{mp}(X,\mathbb{R})$ the max-plus dual of $C(X,\mathbb{R})$, which is the set of all max-plus linear functionals $m: C(X,\mathbb{R}) \to \mathbb{R}$. An element $m \in C_*^{mp}(X,\mathbb{R})$ is called a \emph{Maslov measure} or an \emph{idempotent measure} on $X$.

\begin{definition}\label{def: Idempotent probability as a functional}
	We define \(I(X)\) as the subset of \(C_*^{mp}(X,\mathbb{R})\) consisting of all Maslov measures satisfying \(m(0)=0\). An element \(m \in I(X)\) is called a \emph{Maslov probability} or an \emph{idempotent probability} on \(X\).
\end{definition}

Using the notation \(\mathbb{1}=0\), a Maslov probability satisfies \(m(\mathbb{1})=\mathbb{1}\). Note that for any \(c \in \mathbb{R}\) and \(\mu \in I(X)\), we have \(\mu(c)=c\). Another consequence of this definition is that an idempotent probability is an order-preserving functional. That is, if \(\varphi \leq \psi\), then \(\mu(\varphi) \leq \mu(\psi)\) for \(\varphi, \psi \in C(X, \mathbb{R})\).

Idempotent measures are well known to be closely connected with upper semi-continuous (u.s.c.) functions. The support of a u.s.c. function \(\lambda: X \to \mathbb{R}_{\max}\) is defined as the closed set
\[
\supp(\lambda) := \{ x \in X \mid \lambda(x) \neq -\infty \}.
\]
Hereafter, \( U(X, \mathbb{R}_{\max}) \) denotes the set of u.s.c. functions \(\lambda: X \to \mathbb{R}_{\max}\) such that \( \supp(\lambda) \neq \varnothing \). Additionally, we use \( \bigoplus_{x \in X} \) to denote \( \sup_{x \in X} \). Since \( X \) is compact, every function \( \lambda \in U(X, \mathbb{R}_{\max}) \) achieves its supremum.

The following result is well-established in the literature. For compact metric spaces, refer to \cite{MO24}, and for a different context, see \cite{KM89}.

\begin{theorem}\label{teo : densidade Maslov} A functional $\mu:C(X, \mathbb{R}) \to  \mathbb{R} $ is an idempotent measure if and only if there exists $\lambda\in U(X, \mathbb{R}_{\max})$ satisfying
	\[\mu(\psi) = \bigoplus_{x\in X} \lambda(x)\odot \psi(x),\,\,\,\,\forall \,\psi \in C(X,\mathbb{R}).\]
There is a unique such function $\lambda$ in $U(X, \mathbb{R}_{\max})$  and  $\mu\in I(x)$ if and only if $\oplus_{x\in X}\lambda(x) = 0.$
\end{theorem}

The unique upper semi-continuous function \(\lambda\) presented in Theorem \ref{teo : densidade Maslov} is called the \textit{density} of \(\mu\). We use the notation \(\mu = \bigoplus_{x \in X} \lambda(x) \odot \delta_x \in C^*(X, \mathbb{R})\), where \(\delta_x(\psi) = \psi(x)\). The \emph{support} of an idempotent measure \(\mu = \bigoplus_{x \in X} \lambda(x) \odot \delta_x\) in \(X\) refers to the support of its density.

Following canonical references (see \cite{Rep10}, \cite{Zar}, and \cite{ZAI}), we endow \(I(X)\) with the topology of pointwise convergence \(\tau_p\), whose subbase consists of sets \(V(\mu, f, \varepsilon) := \{\nu \in I(X) : |\nu(f) - \mu(f)| < \varepsilon\}\), where \(\mu \in I(X)\), \(f \in C(X)\), and \(\varepsilon > 0\). Note that \(I(X)\) is compact with respect to \(\tau_p\) provided \(X\) is compact (see, for example, \cite[Theorem 5.3]{Rep10}). Clearly, \(\mu_n \to \mu\) with respect to \(\tau_p\) if \(\mu_n(f) \to \mu(f)\) for all \(f \in C(X)\).

\begin{definition}\label{def:ucCIFS}
Let \((X, d_X)\) be a compact metric space. A \emph{uniformly contractible countable iterated function system} (ucCIFS) \(R = (X, (\phi_j)_{j \in \mathbb{N}})\) is a countable family of maps \(\{\phi_j : X \to X \mid j \in \mathbb{N}\}\) satisfying: there exists \(0 < \gamma < 1\) such that
\begin{equation}\label{eq: gamma contraction}
	d_X(\phi_j(x_1), \phi_j(x_2)) \leq \gamma \cdot d_X(x_1, x_2), \quad \forall j \in \mathbb{N}, \forall x_1, x_2 \in X.
\end{equation}
\end{definition}

From \eqref{eq: gamma contraction}, each map \(\phi_j\) is uniformly \(\gamma\)-Lipschitz continuous, implying that \(\operatorname{Lip}(R) \leq \gamma\).

We recall from Secelean's work \cite{Sec01a} that a countable IFS \(R = (X, (\phi_j)_{j \in \mathbb{N}})\), formed by a sequence of uniform contraction maps, has a unique attractor. This attractor is a compact set \(A\) such that \(\overline{\bigcup_{j \in \mathbb{N}} \phi_j(A)} = A\). Secelean also showed that this attractor \(A\) can be approximated by the attractors \(A_n\) of the partial systems \(R_n = (X, (\phi_j)_{j \in \{1, \ldots, n\}})\), which are formed by considering the first \(n\) maps. Specifically, \(\bigcup_{j \in \{1, \ldots, n\}} \phi_j(A_n) = A_n\) and \(A_n \to A\) with respect to the Hausdorff metric.

We will demonstrate that a similar result holds for a countable max-plus IFS with constant weights.

\begin{definition} \label{def:mpIFS with measures}
Let \((X, d_X)\) be a compact metric space. A \emph{max-plus countable iterated function system} (mpCIFS) \(S = (X, (\phi_j)_{j \in \mathbb{N}}, (q_j)_{j \in \mathbb{N}})\) is a uniformly contractive CIFS \((X, (\phi_j)_{j \in \mathbb{N}})\) equipped with a normalized family of weights \((q_j)_{j \in \mathbb{N}}\), where \(\{q_j \leq 0 \mid j \in \mathbb{N}\}\) satisfies:
\begin{equation}\label{eq: q normalized}
	\bigoplus_{j \in \mathbb{N}} q_j = 0.
\end{equation}
\end{definition}

We also use the compact notation \(S = (X, \phi, q)\) to denote a mpCIFS.

\begin{example}\label{ex:example mpCIFS basic}
   Consider \(X = [0, 1]\) with the usual metric and \(\phi_j(x) = \frac{1}{j+2} x + \frac{1}{j+2}\). Here, \(\gamma = \frac{1}{3}\) is a common contraction rate because \(\sup \operatorname{Lip}(\phi_j) = \sup \frac{1}{j+2} = \frac{1}{3}\). Define the weights by \(q_j = -\frac{1}{j+2}\). Then, \(S = (X, \phi, q)\) is clearly a mpCIFS because \(\bigoplus_{j \in \mathbb{N}} q_j = \sup_{j \in \mathbb{N}} -\frac{1}{j+2} = 0\).
\end{example}

\begin{definition}\label{def: functor}
   For each \(j \in \mathbb{N}\) and \(\mu \in I(X)\), we define the operation \(I_j(\mu)\) for any \(f \in C(X, \mathbb{R})\) as follows:
   \[
   I_j(\mu)(f) := q_j \odot \mu(f \circ \phi_j).
   \]
\end{definition}

\begin{lemma}\label{lem:partial transfer operator}
   To each mpCIFS, $S=(X, \phi, q)$, $\mu \in I(X)$   and $f \in C(X, \mathbb{R})$  consider the functional
   $$\Lambda_{\mu}(f):=\lim_{n \to \infty} \bigoplus_{j=1}^{n} I_{j}(\mu)(f).$$
   Then, $\Lambda_{\mu}= \bigoplus_{j \in \mathbb{N}} I_{j}(\mu) \in I(X)$ is well defined.
\end{lemma}
\begin{proof}
    Consider the sequence
   $$n \to a_n:=\bigoplus_{j=1}^{n} I_{j}(\mu)(f) \in \mathbb{R}$$
   then it is increasing and bounded from above by $\mu(|f|)$. Indeed,
   $$a_{n+1}=\bigoplus_{j=1}^{n+1} I_{j}(\mu)(f)=\max\left\{\bigoplus_{j=1}^{n} I_{j}(\mu)(f), \; I_{n+1}(\mu)(f)\right\} \geq  \bigoplus_{j=1}^{n} I_{j}(\mu)(f)=a_{n}$$
   and, for each $j$ we obtain
   $$I_{j}(\mu)(f)= q_{j} \odot \mu(f\circ \phi_{j}) \leq \mu(|f|)+ (\bigoplus_{j\in \mathbb{N}} q_j) = \mu(|f|)$$
   so $a_{n} \leq  \mu(|f|)$. Thus the limit does exits and equals the supremum. Define $\displaystyle\Lambda_{\mu}(f):=\lim_{n \to \infty} a_{n} = \bigoplus_{j \in \mathbb{N}} I_{j}(\mu)(f)$, then the functional is well defined. In particular, $\bigoplus_{j=1}^{n} I_{j}(\mu) \stackrel{\tau_p}{\to} \Lambda_{\mu}$.

   To show that $\Lambda_{\mu} \in I(X)$ we can easily check the three properties:
   \begin{enumerate}
        \item $\displaystyle\Lambda_{\mu}(0)= \lim_{n \to \infty} \bigoplus_{j=1}^{n} I_{j}(\mu)(0)= 0$, because $\mu(0)= 0$;
		\item
$$\Lambda_{\mu}(c \odot f)= \lim_{n \to \infty} \bigoplus_{j=1}^{n} I_{j}(\mu)(c \odot f) = \lim_{n \to \infty} \bigoplus_{j=1}^{n} \mu(q_j \odot c \odot f(\phi_j)) =$$
$$=\lim_{n \to \infty} c \odot \bigoplus_{j=1}^{n} \mu(q_j \odot  f(\phi_j))=  c \odot \Lambda_{\mu}( f),$$
$\forall c \in \mathbb{R}$ and $\forall f \in C(X,\mathbb{R})$;
		\item
$$\Lambda_{\mu}(f \oplus g)= \lim_{n \to \infty} \bigoplus_{j=1}^{n} I_{j}(\mu)(f \oplus g)= \lim_{n \to \infty} \bigoplus_{j=1}^{n} \mu(q_j \odot [f(\phi_j) \oplus g(\phi_j)]) =$$
$$= \lim_{n \to \infty} \bigoplus_{j=1}^{n} [\mu(q_j \odot f(\phi_j) ) \oplus \mu(q_j \odot g(\phi_j))] = \Lambda_{\mu}( f) \oplus \Lambda_{\mu}( g),$$
$\forall f, g \in C(X,\mathbb{R})$.
	\end{enumerate}
\end{proof}

\begin{definition}\label{def:Markov operator}
	To each mpCIFS $S=(X, \phi, q)$, we assign the operators $M_{\phi,q,n}:C_*^{mp}(X,\mathbb{R}) \to C_*^{mp}(X,\mathbb{R})$, and $M_{\phi,q,n}: I(X) \to C_*^{mp}(X,\mathbb{R})$ defined by
\[M_{\phi,q,n}(\mu)(f):= \bigoplus_{j=1}^{n} I_{j}(\mu)(f )
\]
and
	\begin{equation}\label{eq:partial Markov operator}
		M_{\phi,q}(\mu)(f):= \Lambda_{\mu}(f)=\lim_{n \to \infty}M_{\phi,q,n}(\mu)(f) = \bigoplus_{j \in J} I_{j}(\mu)(f),
	\end{equation}
	for any $f \in C(X, \mathbb{R})$. An idempotent probability  $\mu \in I(X)$ is called invariant (with respect to the mpCIFS) if $M_{\phi,q}(\mu)=\mu$.
\end{definition}

We observe that \( M_{\phi,q,n} \) is not a normalized max-plus operator unless \(\bigoplus_{j=1}^{n} q_{j} = 0\) (i.e., if \( q_{j} = 0 \) for some \( j \leq n \)). However, we can normalize it by subtracting \(\alpha_n := \bigoplus_{j=1}^{n} q_{j} \leq 0\) from each \( q_{j} \). Then, the operator \( M_{\phi,q,n} \) can be replaced by \( M_{\phi,\tilde{q},n}: I(X) \to I(X) \) defined by
\[
M_{\phi,\tilde{q},n}(\mu)(f) := \bigoplus_{j=1}^{n} \tilde{q}_{j} \odot \mu(f \circ \phi_{j}),
\]
where \(\tilde{q}_{j} := q_{j} - \alpha_{n}\). Thus, \( M_{\phi,q,n}(\mu) = \alpha_{n} \odot M_{\phi,\tilde{q},n}(\mu) \).

From \cite{MZ} and \cite{dCOS20fuzzy,ExistInvIdempotent}, we know that there exists a unique idempotent measure \( \mu_{n} \in I(X) \) such that \( M_{\phi,\tilde{q},n}(\mu_{n}) = \mu_{n} \) or \( M_{\phi,q,n}(\mu_{n}) = \alpha_{n} \odot M_{\phi,\tilde{q},n}(\mu_{n}) = \alpha_{n} \odot \mu_{n} \). It is also easy to see that \(\alpha_n\) is increasing and \(\displaystyle \lim_{n \to \infty} \alpha_n = 0\). Motivated by this, we define:
\begin{definition} \label{def:partial mpIFS with measures}
Let $(X,d_X)$  be a compact metric space and $S=(X, (\phi_j)_{j\in \mathbb{N}}, (q_j)_{j\in \mathbb{N}})$ a uniformly contractive mpCIFS. For each $n \in \mathbb{N}$ we define the partial  mpIFS induced by $S$ as a (finite) normalized mpIFS $S_{n}=(X, (\phi_j)_{j \leq n}, (\tilde{q}_j)_{j\leq n})$ where $\alpha_{n}:= \bigoplus_{j=1}^{n} q_{j}$ and
$\tilde{q}_{j}:=q_{j}- \alpha_{n}$.
\end{definition}
As a consequence of Definition~\ref{def:ucCIFS}, we know that the maps $\phi_j, \; 1 \leq j \leq n$, are uniformly $\gamma$-Lipschitz continuous, with respect to $j$ thus $\operatorname{Lip}(S_{n}) \leq \gamma$.

From \cite{MO24} we know that for each mpIFS $S_{n}=(X, \phi, \tilde{q})$ we assign the following operators:\newline
	1.  $\mathcal{L}_{\phi,\tilde{q},n}: C(X, \mathbb{R})   \to C(X, \mathbb{R})  $, defined by
	\begin{equation}\label{eq:mp ruelle operator}
		\mathcal{L}_{\phi,\tilde{q},n}(f)(x):=\bigoplus_{1\leq j\leq n} \tilde{q}_{j}(x)\odot f(\phi_{j}(x)).
	\end{equation}
	2.  $M_{\phi,\tilde{q},n}: I(X) \to I(X),$ defined by
	\begin{equation}\label{eq:Markov operator}
		M_{\phi,\tilde{q},n}(\mu):= \bigoplus_{1\leq j\leq n}\mu( \tilde{q}_{j}\odot (f\circ \phi_{j})).
	\end{equation}
	3.  $L_{\phi,\tilde{q},n}: U(X, \mathbb{R}_{\max}) \to U(X, \mathbb{R}_{\max}),$ defined by
	\begin{equation}\label{eq:transfer operator}
		L_{\phi,\tilde{q},n}(\lambda)(x):=\bigoplus_{(j,y)\in \phi^{-1}(x)} \tilde{q}_{j} \odot  \lambda(y).
	\end{equation}

Considering these operators and an idempotent probability   $\mu=\bigoplus_{x\in X}\lambda(x)\odot\delta_x\in I(X)$, we have that $M_{\phi,\tilde{q},n}(\mu)= \bigoplus_{x \in X} L_{\phi,\tilde{q},n} (\lambda)(x) \odot \delta_x$, that is, $M_{\phi,\tilde{q},n}(\mu)$ has density $L_{\phi,\tilde{q},n} (\lambda)$ where $\lambda$ is the density of $\mu$. Furthermore $M_{\phi,\tilde{q},n}(\mu)(f) = \mu(\mathcal{L}_{\phi,\tilde{q},n}(f)),$ for any $f \in C(X, \mathbb{R})$, that is, $M_{\phi,\tilde{q},n}$ is the max-plus dual of $\mathcal{L}_{\phi,\tilde{q},n}$.

\section{The existence of the attractor by contractivity and approximation by partial attractors}\label{sec:approx by contract}

There are several ways to introduce topologies on \( I(X) \) other than the pointwise topology, \(\tau_p\).

We recall from \cite{ExistInvIdempotent} that for each \( a > 0 \) and \(\mu, \nu \in I(X)\), the pseudometric \(d_a\) is defined as
\[
d_a(\mu, \nu) = \sup \{ |\mu(g) - \nu(g)| : g \in \operatorname{Lip}_a(X) \},
\]
where \(\operatorname{Lip}_a(X)\) is the family of maps \( g: X \to \mathbb{R} \) with \(\operatorname{Lip}(g) \leq a\). In \cite[Theorem 4.1]{Rep10}, it is proven that \( d_a \) are continuous pseudometrics for each \( a \in \mathbb{N} \), and that \(\tilde{d}\) defined by
\[
\tilde{d}(\mu, \nu) := \sum_{i=1}^\infty \frac{d_i(\mu, \nu)}{i \cdot 2^i}, \quad \mu, \nu \in I(X)
\]
is a metric on \( I(X) \) that generates the canonical topology \(\tau_p\).

In \cite{ExistInvIdempotent}, the authors show that the idempotent Markov operator for Banach contractive normalized mpIFS is a Banach contraction with respect to the following natural modification of \(\tilde{d}\):

\begin{definition}\label{def: new metric I(X) repovs Filip}
For \(\beta, \tau \in (0,1)\), define \(\tilde{d}_{\beta,\tau}\) by
\[
\tilde{d}_{\beta,\tau}(\mu, \nu) := \sum_{i \in \mathbb{Z}} \frac{\tau^{|i|}}{\beta^i} d_{\beta^i}(\mu, \nu),
\]
for any \(\mu, \nu \in I(X)\).
\end{definition}

In a similar manner to \cite{Rep10}, one can show that \(\tilde{d}_{\beta,\tau}\) is a metric that generates the topology \(\tau_p\). We just need to observe that \(\tilde{d}_{\beta,\tau}\) is well-defined and that the family \( d_{\beta^i} \), \( i \in \mathbb{Z} \), is a family of continuous pseudometrics that separates points. In particular, \((X, d_{\beta^i})\) is compact because it generates the compact topology \(\tau_p\), so \((X, d_{\beta^i})\) is complete.

The following theorem from \cite{ExistInvIdempotent} provides an alternative proof of the existence of an invariant idempotent measure for Banach contractive normalized mpIFS.

\begin{theorem}\cite[Theorem 4.1]{ExistInvIdempotent}\label{thm: contractivity mpIFS Filip}
Assume that $S_{n}=(X,(\phi_j)_{j=1}^L,(\tilde{q}_j)_{j=1}^n)$ is a Banach contractive normalized mpIFS. Let
$$
\gamma:=\max\{\operatorname{Lip}(\phi_j):j=1,...,n\}=\operatorname{Lip}(S_{n})
$$
and choose $\tau \in(\gamma,1)$. Then $M_{\phi,\tilde{q}, n}$ is Banach contraction with respect to $\tilde{d}_{\gamma,\tau}$. More precisely,
$$
\operatorname{Lip}(M_{\phi,\tilde{q},n})\leq \frac{\gamma}{\tau}<1.
$$
\end{theorem}

\begin{theorem}\label{thm:contraction mpCIFS}
   Let $S=(X, (\phi_j)_{j\in \mathbb{N}}, (q_j)_{j\in \mathbb{N}})$ be a $\gamma$-uniformly contractive mpCIFS and $\tau \in(\gamma,1)$. Then
   $$
     \operatorname{Lip}(M_{\phi,q})\leq \frac{\gamma}{\tau}<1,
   $$
   that is, $M_{\phi,q}$ is a Banach contraction (with respect to the metric $\tilde{d}_{\gamma,\tau}$), so there exists a unique $\nu \in I(X)$ such that $M_{\phi,q}(\nu)=\nu$ and the sequence $M_{\phi,q}^{n}(\mu) \to \nu$ w.r.t. the metric $ \tilde{d}_{\gamma,\tau}$, for any initial point $\mu \in I(X)$.
\end{theorem}
\begin{proof}
   Recall that
   $$M_{\phi,q}(\mu)(g)=\lim_{n \to \infty} M_{\phi,q,n}(\mu)(g)= \lim_{n \to \infty} \alpha_{n} \odot M_{\phi,\tilde{q},n}(\mu)(g) $$
   for any $g \in C(X,\mathbb{R})$, and $\alpha_n=\bigoplus_{j=1}^{n}q_j$ is increasing with $\displaystyle\lim_{n\to \infty} \alpha_n=0$.

   Fixed $\varepsilon >0$ and $g\in\operatorname{Lip}_a(X)$ there exists $N_0 \in \mathbb{N}$ such that for any $n \geq N_0$ we get
   $$M_{\phi,q}(\mu)(g) -\varepsilon  \leq  \alpha_{n} \odot M_{\phi,\tilde{q},n}(\mu)(g) \leq M_{\phi,q}(\mu)(g)$$
   and
   $$M_{\phi,q}(\nu)(g) -\varepsilon  \leq  \alpha_{n} \odot M_{\phi,\tilde{q},n}(\nu)(g) \leq M_{\phi,q}(\nu)(g).$$
   Subtracting this two equations and denoting $\Delta:= M_{\phi,q}(\mu)(g) - M_{\phi,q}(\nu)(g)$ and $\Delta_n:= M_{\phi,\tilde{q},n}(\mu)(g) - M_{\phi,\tilde{q},n}(\nu)(g)$, we obtain
   $$|\Delta - \Delta_n| < \varepsilon.$$
   Now, take $g\in\operatorname{Lip}_a(X)$ such that
   $$d_a(M_{\phi,q}(\mu),M_{\phi,q}(\nu))-\varepsilon< |M_{\phi,q}(\mu)(g)-M_{\phi,q}(\nu)(g)|= |\Delta|\leq$$
   $$\leq |\Delta - \Delta_n|+| \Delta_n|< \varepsilon +| \Delta_n| \leq  \varepsilon + d_a(M_{\phi,\tilde{q},n}(\mu),M_{\phi,\tilde{q},n}(\nu)),$$
   so that
   $$d_a(M_{\phi,q}(\mu),M_{\phi,q}(\nu)) \leq 2\varepsilon + d_a(M_{\phi,\tilde{q},n}(\mu), M_{\phi,\tilde{q},n}(\nu)).$$
   In the next, we take $a:= \gamma^i $ so that $n:=n_i \geq N_0$ is fixed.
   We recall that,
   $$
\tilde{d}_{\gamma,\tau}(M_{\phi,q}(\mu),M_{\phi,q}(\nu))=
\sum_{i\in\mathbb{Z}}\frac{\tau^{|i|}}{\beta^i}d_{\gamma^i}(M_{\phi,q}(\mu),M_{\phi,q}(\nu))\leq
$$
$$\leq \sum_{i \in\mathbb{Z}}\frac{\tau^{|i|}}{\gamma^i} \left[ 2\varepsilon + d_{\gamma^i}(M_{\phi,\tilde{q},n_i}(\mu),M_{\phi,\tilde{q},n_i}(\nu)) \right].$$

Note that $\operatorname{Lip}(S_{n_i}) \leq \gamma$ uniformly with respect to $n_i$ and by Theorem~\ref{thm: contractivity mpIFS Filip} we obtain
$$d_{\gamma^i}(M_{\phi,\tilde{q},n_i}(\mu),M_{\phi,\tilde{q},n_i}(\nu))\leq d_{\gamma^{i+1} }(\mu,\nu).$$
Using this equation we obtain,
$$
\tilde{d}_{\gamma,\tau}(M_{\phi,q}(\mu),M_{\phi,q}(\nu)) \leq
 2\varepsilon \sum_{i \in\mathbb{Z}}\frac{\tau^{|i|}}{\gamma^i} +
\sum_{i \in\mathbb{Z}}\frac{\tau^{|i|}}{\gamma^i}  d_{\gamma^{i+1} }(\mu,\nu)=$$
$$=2\varepsilon \sum_{i \in\mathbb{Z}}\frac{\tau^{|i|}}{\gamma^i} + \frac{\gamma}{\tau} \,
\sum_{r \in\mathbb{Z}}\frac{\tau^{|r|}}{\gamma^{r}}  d_{\gamma^{r} }(\mu,\nu)= 2\varepsilon \sum_{i \in\mathbb{Z}}\frac{\tau^{|i|}}{\gamma^i} + \frac{\gamma}{\tau} \, \tilde{d}_{\gamma,\tau}(\mu,\nu).$$
Since $\varepsilon>0$ is arbitrary we proved that
$$
\tilde{d}_{\gamma,\tau}(M_{\phi,q}(\mu),M_{\phi,q}(\nu)) \leq  \frac{\gamma}{\tau} \, \tilde{d}_{\gamma,\tau}(\mu,\nu).$$
\end{proof}

The next theorem shows that the attractor of a mpCIFS can be pointwisely approximated by the attractors of partial mpIFSs, that is, the unique idempotent probabilities $\mu_{n} \in I(X)$ satisfying $M_{\phi,\tilde{q},n}(\mu_{n})=\mu_{n}$.

 \begin{theorem}\label{thm:contraction mpCIFS approximation}
   Let $S=(X, (\phi_j)_{j\in \mathbb{N}}, (q_j)_{j\in \mathbb{N}})$ be a $\gamma$-uniformly contractive mpCIFS and  $\nu \in I(X)$ the unique idempotent probability  such that $M_{\phi,q}(\nu)=\nu$ (given by Theorem~\ref{thm:contraction mpCIFS}). Consider $\mu_{n} \in I(X)$ the unique idempotent probability such that $M_{\phi,\tilde{q},n}(\mu_{n})=\mu_{n}$. Let  $\mu_{n_i}$
   for $i \in \mathbb{N}$, be any convergent subsequence (w.r.t. the $\tau_p$ topology). Then, $\mu_{n_i} \stackrel{\tau_p}{\to} \nu$, in particular $\mu_{n}$ is convergente w.r.t. the $\tau_p$ topology.
\end{theorem}
\begin{proof}
   Fix $g \in C(X,\mathbb{R})$ and suppose that $\mu_{n_i} \stackrel{\tau_p}{\to} \mu$. Since $\displaystyle M_{\phi,q}(\nu)(g)=\lim_{n \to \infty } \alpha_{n} \odot M_{\phi,\tilde{q},n}(\mu)(g)$ we can find, for any $\varepsilon>0$ a natural number $N_{\varepsilon} \in \mathbb{N}$ such that
   $$| M_{\phi,q}(\mu)(g) - \alpha_{n} \odot M_{\phi,\tilde{q},n}(\mu)(g) |< \frac{\varepsilon}{5},$$
   for any $n \geq N_{\varepsilon}$.

   On the other hand $\alpha_n:=\bigoplus_{j=1}^{n} q_{j} \leq 0 $ is increasing and $\lim_{n\to \infty} \alpha_n=0$. So we can assume
   $-\alpha_n < \frac{\varepsilon}{5},$  for any $n \geq N_{\varepsilon}$.

   Also, from the convergence $\mu_{n_i} \stackrel{\tau_p}{\to} \mu$, we obtain
   $| \mu_{n_i}(g) - \mu(g) |< \frac{\varepsilon}{5},$
   for any $n \geq N_{\varepsilon}$.

   Finally, from the continuity of $M_{\phi,\tilde{q},n}$ with respect to $\mu$, we have
   $$| M_{\phi,\tilde{q},n}(\mu)(g) - M_{\phi,\tilde{q},n}(\mu_{n})(g) |< \frac{\varepsilon}{5},$$
   for any $n \geq N_{\varepsilon}$.

   Thus, we can write, for $n=n_i$,
   $$| M_{\phi,q}(\mu)(g) - \mu(g) | = $$
   $$=|M_{\phi,q}(\mu)(g) - \alpha_{n} \odot M_{\phi,\tilde{q},n}(\mu)(g) + \alpha_{n} \odot M_{\phi,\tilde{q},n}(\mu)(g) -  M_{\phi,\tilde{q},n}(\mu_{n})(g)+ M_{\phi,\tilde{q},n}(\mu_{n})(g) - \mu(g) | \leq $$
   $$\leq |M_{\phi,q}(\mu)(g) - \alpha_{n} \odot M_{\phi,\tilde{q},n}(\mu)(g)| +  |\alpha_{n}| +$$ $$+ |M_{\phi,\tilde{q},n}(\mu)(g) -  M_{\phi,\tilde{q},n}(\mu_{n})(g)|+ |M_{\phi,\tilde{q},n}(\mu_{n})(g) - \mu(g)| < \frac{4}{5}\varepsilon < \varepsilon,$$
   because $M_{\phi,\tilde{q},n}(\mu_{n})=\mu_{n}$.
    Since $\varepsilon>0$ is arbitrary we obtain $| M_{\phi,q}(\mu)(g) - \mu(g) | = 0$ for $g \in C(X,\mathbb{R})$ , that is,  $ M_{\phi,q}(\mu)= \mu$. As $\nu$ is the unique solution of this equation we conclude that $\mu_{n_i} \stackrel{\tau_p}{\to} \nu$, for any subsequence, and so $\mu_{n} \stackrel{\tau_p}{\to} \nu$.
\end{proof}

\begin{remark}
An open question is if we can estimate the convergence rate $\mu_{n} \to \nu$ with respect to some metric on $I(X)$ such as the ones used in \cite{dCOS20fuzzy} ($d_\theta$) or \cite{ExistInvIdempotent} ($\tilde{d}_{\gamma,\tau}$). For instance, if one can show that
$$
\tilde{d}_{\gamma,\tau}(M_{\phi,q}(\mu_{n}),\mu_{n+1}) \leq \varepsilon_{n},$$
and $\varepsilon_{n} \to 0$, then, Ostrowski's stability theorem for Banach contractions, \cite{Ostr1967}, gives a convergence rate, since $M_{\phi,q}$ is a contraction with respect to $\tilde{d}_{\gamma,\tau}$ by Theorem~\ref{thm:contraction mpCIFS}.
\end{remark}

\section{Fuzzy approach for idempotent invariant measures}\label{sec: fuzzy approach}

From \cite{dCOS20fuzzy}, we know that the scale map \(\theta: [-\infty, 0] \to [0,1]\) defined by \(\theta(t) = e^t\) (among other possible choices) induces a bijection \(\Theta\) between \(I(X)\) and \(\mathcal{F}_{X}^*\) (the set of normal upper semicontinuous fuzzy sets), given by
\[
u(x) = \Theta(\mu)(x) = \theta(\lambda_{\mu}(x)),
\]
for any \(\mu \in I(X)\) with density \(\lambda_{\mu}\).

Moreover, consider a mpIFS \(S = (X, (\phi_j)_{j \in J}, (q_j)_{j \in J})\), where \(J = \{1, \ldots, n\}\), \(\phi_{j}: X \to X\), and \(q_{j} \in [-\infty, 0]\) such that \(\bigoplus_{j \in J} q_{j} = 0\). From \cite{dCOS20fuzzy}, we know that given the associated fuzzy iterated function system (IFZS), \(\mathcal{Z}_S = (X, (\phi_j)_{j \in J}, (\rho_j)_{j \in J})\), where the gray level functions are
\[
\rho_j(t) := \theta(q_j + \theta^{-1}(t)) = e^{q_j} \, t, \quad t \in [0,1],
\]
we have
\[
\Theta \circ M_{\phi,q} = F_{\mathcal{Z}_S} \circ \Theta.
\]
In particular, \(M_{\phi,q}(\mu) = \mu\) if and only if \(F_{\mathcal{Z}_S}(u) = u\) for \(u = \Theta(\mu)\), where the fuzzy fractal operator is
\[
F_{\mathcal{Z}_S}(v) = \bigvee_{1 \leq j \leq n} \rho_j(\phi_j(v)),
\]
for any \(v \in \mathcal{F}_{X}^*\) (see \cite{Oli17} for details on IFZS).

Then, we can introduce a metric \(d_{\theta}\) induced by \(\Theta\) from the metric space of fuzzy sets \((\mathcal{F}_{X}^*, d_{\text{f}})\) (where \(d_{\text{f}}\) is the supremum of the Hausdorff distance between level sets), given by \(d_{\theta}(\mu, \nu) = d_{\text{f}}(\Theta(\mu), \Theta(\nu))\), in such a way that the spaces \((I(X), d_\theta)\) and \((\mathcal{F}_{X}^*, d_{\text{f}})\) are homeomorphic. By \cite[Theorem 3.5]{dCOS20fuzzy}, the space \((I(X), d_\theta)\) is complete since \((X, d)\) is compact. The following two results hold:

\begin{lemma}\cite[Lemma 5.7]{dCOS20fuzzy}\label{filipmetric}
For every  $\mu=\bigoplus_{x\in X}\lambda(x)\odot\delta_x,\;\nu=\bigoplus_{x\in X}\eta(x)\odot\delta_x\in I(X)$, we have
$$
d_\theta(\mu,\nu)=\bigoplus_{\beta\in(-\infty,0]}h(\{x\in X:\lambda(x)\geq \beta\},\{x\in X:\eta(x)\geq \beta\}),
$$
where $h$ is the Hausdorff distance.
\end{lemma}

\begin{proposition}\cite[Proposition 5.8]{dCOS20fuzzy}\label{fil11}
The metric space $(I(X),\, d_\theta)$ is complete and the topology $\tau_\theta$ induced by $d_\theta$ is finer than the  pointwise convergence topology $\tau_p$. In other words, $\tau_p\subset\tau_\theta$.
\end{proposition}

We end this section recalling  an important result from \cite{MO24} who gives a representation for the fuzzy attractor in correspondence with the invariant idempotent probability of a mpIFS (see \cite{dCOS20fuzzy} and  \cite[Theorem 4.7]{MO24}).

\begin{proposition}\label{prop:appl fuzzy}
    Consider $S=(X, (\phi_j)_{j\in J}, (q_j)_{j\in J})$, $J=\{1,\ldots,n\}$, a Banach contractive normalized mpIFS, with constant weights, and the associated IFZS
    $$\mathcal{Z}_S=(X,(\phi_j)_{j\in J},(\rho_j)_{j\in J}),$$
    where $\rho_j(t)=e^{q_j}t,\, t \in [0,1]$ and $\theta=e^t$ is a scale function. Then, the fuzzy attractor of $\mathcal{Z}_S$ satisfies
    $$u(x)= \bigoplus_{(j_1,j_2,j_3,...) \in \pi^{-1}_{n}(x)}e^{(q_1+q_2+q_3+...)},$$
where $\pi_{n}: \Omega_{n} \to X$ is given by
$\displaystyle\pi_{n}(j) = \lim_{k\to+\infty} \phi_{j_1}\circ...\circ\phi_{j_k}(X)$
for any sequence $j=(j_1,j_2, \ldots ) \in \Omega_{n}:=\{1,\ldots,n\}^{\mathbb{N}}$
\end{proposition}
We may consider to apply this result to $S=S_n$ to study the fixed points $M_{\phi,\tilde{q},n}(\mu_{n})=\mu_{n}$. Although, despite the fact that this solution is explicit, we do not have a good numerical procedure to compute it.

\begin{example}\label{ex:example mpCIFS explicit}
    Let us consider $X=[0,1]$ endowed with the usual Euclidian induced metric. For each $j=1,2,\ldots$ we define the map $\phi_j: X \to X $ by
    $$\phi_j(x):= \frac{1+ x}{2^{j}}.$$
    We also consider the weights    $q_j(x):= -\frac{1}{2^j}.$
    Let $S=(X,\phi,q)$ be the associated contractive mpCIFS, which is normalized because $\bigoplus_{j \in \mathbb{N}} -\frac{1}{2^j} =0$.
\begin{figure}[H]
  \centering
  \includegraphics[width=8cm]{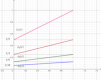}
  \caption{Plot of the maps for the IFS $S_4$.}\label{fig:example1-1}
\end{figure}

We aim to describe its invariant measure  through the invariant measures of the partial induced normalized  mpIFS
$$S_{n}=(X,\phi_{j},\tilde{q}^{n}_{j}),\;  j=1,\ldots,n,$$
where
$\tilde{q}^{n}_{j}= q_{j} -\alpha_{n}$ and $\alpha_{n}= \bigoplus_{1 \leq j \leq n} -\frac{1}{2^j} = -\frac{1}{2^n}$. Let $\mu_{n}$ be the unique invariant idempotent probabilities for $S_{n}$.

Consider the associated IFZS  $\mathcal{Z}_{n}:=\mathcal{Z}_{S_{n}}=(X,(\phi_j)_{j\in J},(\rho_j)_{j\in J}),$  where $\rho_j(t)=e^{q_j}t$, then the fuzzy attractor of $\mathcal{Z}_{S_{n}}$ (see \cite{Oli17} for details on the existence and uniqueness) satisfies
    $$u(x)=F_{\mathcal{Z}_{n}}(u)(x)= \bigvee_{1 \leq j \leq n}\rho_j(\phi_j(u)),$$
    where the Zadeh extension principle (see \cite{Zad} for details on fuzzy sets, also \cite{Oli17} for the general case) holds
    $$\phi_j(u)(x)= \sup_{\phi_j(y)=x} u(y) =
    \left\{
      \begin{array}{ll}
        0, &  x < \frac{1}{2^{j}} \text{ or } x > \frac{1}{2^{j-1}}\\
        u(0), & x=\frac{1}{2^{j}} \\
        u(1), & x=\frac{1}{2^{j-1}}\\
        u(2^{j} x -1), & x  \in (\frac{1}{2^{j}}, \, \frac{1}{2^{j-1}})
      \end{array}
    \right.
    $$
By our previous construction, if we can find the fixed point $u(x)$ satisfying $u(x)=F_{\mathcal{Z}_{n}}(u)(x)$ then the idempotent probability $\mu_n$ with density $\lambda_{\mu_{n}}(x)=\theta^{-1}(u(x)),$ is necessarily invariant for the mpIFS $S_n$.
\end{example}

\section{Higuchi fractal dimension via numerical approximations}\label{Hichi approximations}

The key to approximating the fixed point of a max-plus contractive iterated function system (mpCIFS) \(S\) lies in Theorem~\ref{thm:contraction mpCIFS approximation}. This theorem allows us to approximate the fixed point by using the fixed points \(\mu_n\) of regular max-plus iterated function systems (mpIFS) \(S_n\). By employing the algorithms from \cite{dCOS20fuzzy}, we can generate visual representations of the approximate attractor \(\mu_n\) with a guaranteed approximation rate.

Following \cite{dCOS20fuzzy}, we recall that an \(\varepsilon\)-net in \((X, d)\) is a subspace \((\hat{X}, d)\) equipped with a projection map \(r: X \to \hat{X}\) such that \(d(x, r(x)) < \varepsilon, \, \forall x \in X\). Naturally, each map \(\phi: X \to X\) induces a discrete version \(\hat{\phi}: \hat{X} \to \hat{X}\) given by \(\hat{\phi}(y) = r(\phi(y))\) for any \(y \in \hat{X}\).

Let \(\hat{S}_n\) be the discretization of \(S_n\), producing the discrete operator \(\hat{M}_{\hat{\phi}, \tilde{q}, n}\). Since \(S_n\) is Banach contractive, after \(N\) iterations, \(\hat{M}_{\hat{\phi}, \tilde{q}, n}^N(\eta)\), starting from the initial measure \(\eta\), approximates the actual attractor \(\mu_n\) of \(M_{\phi, \tilde{q}, n}\).

Using the identification of idempotent measures as fuzzy sets, as explained in Section~\ref{sec: fuzzy approach}, one can prove from \cite[Theorem 6.3]{dCOS21} or \cite{dCOS20fuzzy} that there exists a sufficiently large \(N\) to ensure that \(\mu_n\) is approximated with resolution \(\delta > 0\), such that
\[
d_{\theta}(\mu_n, \hat{M}_{\hat{\phi}, \tilde{q}, n}^N(\eta)) < \delta := \frac{2 \varepsilon}{1 - \operatorname{Lip}(S_n)} = \frac{2 \varepsilon}{1 - \gamma}.
\]

The following algorithm, presented in \cite[Section 7]{dCOS20fuzzy}, outlines the steps to achieve this approximation.

{\tt
\begin{tabbing}
aaa\=aaa\=aaa\=aaa\=aaa\=aaa\=aaa\=aaa\= \kill
     \> \texttt{DeterminIFSIdempMeasureDraw}($S$)\\
     \> {\bf input}: \\
     \> \>  \> $K \subseteq {X}$, any finite and {nonempty} subset (a list of points in ${X}$).\\
     \> \>  \> $\nu$, any discrete idempotent measure such that $\supp (\nu)= K$. \\
     \> {\bf output:} \\
     \> \> \>  A bitmap representing an approximation of the attractor.\\
     \> \> \> A bitmap  image representing a discrete invariant idempotent measure,\\ \>\>\>
 with a gray scale color bar indicating\\ \>\>\> the mass of each pixel.\\
     \> {\bf Initialize}  $\displaystyle\mu:={\bigoplus_{x \in X}{-\infty\odot\delta_{x}}} $,  $W:=\varnothing$\\
     \> {\bf for n from 1 to N do}\\
     \>  \>{\bf for $\ell$ from 1 to $\operatorname{Card}({K})$ do}\\
     \>  \>  \>  \>{\bf for $j$ from $1$ to $L$ do}\\
     \>  \>  \>  \>  \>${x}:={K}[\ell]$\\
     \>  \>  \>  \>  \>${y}:={\phi}_{j}({x})$\\
     \>  \>  \>  \>  \>${W:=W\cup\{y\}}$\\
     \>  \>  \>  \>  \> {$\mu({y}):=\max\{\mu({y}) ,  q_{j} + \nu({x})\}$}\\
     \>  \>  \>  \>{\bf end do} \\
     \>  \>{\bf end do} \\
     \> \> $\nu:=\mu$, $\displaystyle\mu:=\bigoplus_{x \in X}-\infty\odot\delta_{x}$, {$K:=W$ and $W:=\varnothing$}\\
     \>{\bf end do}\\
     \>{\bf return: Print} {$K$} and $\nu$ \\
\end{tabbing}}

The partial attractors \(\mu_n = \bigoplus_{x \in X} \lambda_n(x) \odot \delta_x\) are approximated by a finite density \(\hat{\lambda}_n(x)\), which is a series of values in \([- \infty, 0]\). As these values are not suitable for calculating the fractal dimension, we will use the approach in Section~\ref{sec: fuzzy approach} to transform them into a discrete fuzzy set \(\hat{u}_n(x) := e^{\hat{\lambda}_n(x)}\), resulting in a series of values in \([0, 1]\). We will then compute the Higuchi fractal dimension following \cite{Hig88}, as well as its generalizations for two-dimensional data series according to \cite{Spa14} or \cite{ASR15}.

We note that we will not focus on the numerical precision with which we can estimate \(\mu_n\) through \(\hat{\lambda}_n(x)\). This issue is thoroughly addressed in \cite{dCOS20fuzzy} for mpIFSs or in \cite{dCOS21} when identifying mpIFS with its fuzzy counterpart, as explained in Section~\ref{sec: fuzzy approach}. Additionally, since \(\mu_n \stackrel{\tau_p}{\to} \nu\), we do not control the overall speed of convergence. This aspect is left for future investigation. Therefore, our analysis will concentrate on the fractal behavior of \(\hat{\lambda}_n(x)\) using the Higuchi dimension as a measure of the complexity of the global attractor \(\nu\) of our mpCIFS.

\begin{figure}[H]
  \centering
  \includegraphics[width=10cm]{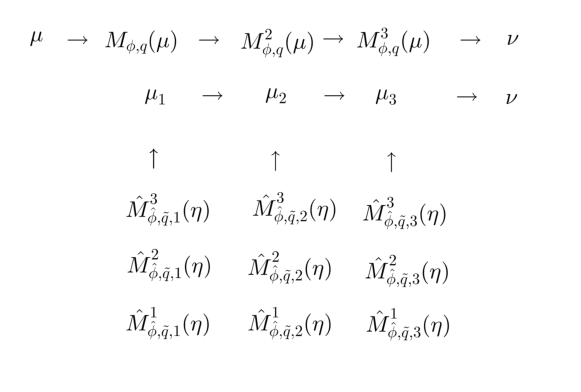}
  \caption{Approximation scheme to obtain $\nu$, where $\mu$ and $\eta$ are arbitrary initial idempotent measures.}\label{fig:approx_scheme}
\end{figure}

\begin{definition}\label{def:eps HFD}
   Let $S$ be a mpCIFS, $n \in \mathbb{N}$ and $\delta, \varepsilon>0$. We define the $(n, \varepsilon, \delta)$-fractal dimension of the attractor of $S$ as the fractal dimension of the discrete fuzzy set $\hat{u}_n(x):=e^{\hat{\lambda}_n(x)}$ where $\hat{\lambda}_n(x)$ is a $\delta$-approximation of the attractor of the discrete mpIFS obtained from $S_n$ by an $\varepsilon$-net $(\hat{X},d)$ in $(X,d)$.
\end{definition}

We will work mostly with $X=[0,1]$ or $X=[0,1]^2$ and uniform $\varepsilon$-nets. In both cases the discrete fuzzy sets will be tim series  of data whose fractal dimension is given by the Higuchi fractal dimension.

\subsection{1-D idempotent fractals and the Higuchi fractal dimension}\label{subsec:higuchi 1-d}
Higuchi's fractal dimension (HFD), introduced by T. Higuchi in 1988 \cite{Hig88}, is a method for determining the fractal characteristics of a time series. Over the years, HFD has found numerous applications across various fields, including the analysis of MGE, MRI, and EEG signals, digital image analysis, pattern recognition, neuroscience, histology, botany, medicine, and physics. Since 2011, significant efforts have been made to generalize HFD for other types of statistical data sets, as seen in works such as \cite{Aha11}, \cite{Spa14}, and \cite{ASR15}, among many others.

We begin by recalling the original Higuchi fractal dimension for completeness. Given a time series \(X: \{1, \ldots, N\} \to \mathbb{R}\), the goal is to estimate the lengths of \(X\) with respect to different scales. The algorithm, detailed in \cite{LM20}, investigates mathematically precise conditions under which the Higuchi method provides the correct value for the box-counting dimension of a function's graph. The procedure, as described in Higuchi's original paper \cite{Hig88}, is as follows:

\begin{figure}[H]
{\tt
\begin{tabbing}
aaa\=aaa\=aaa\=aaa\=aaa\=aaa\=aaa\=aaa\= \kill
    \>  \texttt{1D HFD algorithm} \\
    \> Input: Choose $2\leq N \in \mathbb{N}$, $X:\{1,...,N\} \to \mathbb{R}$ and $2 \leq k_{max} \leq \lceil\frac{N}{2}\rceil$\\
     \> Output: Higuchi Fractal Dimension $HFD(X, N, k_{max})$\\
     \>  {\bf for } $k$ {\bf from } $1$ to $k_{max}$ {\bf do } \\
     \> \>   {\bf for } $m$ {\bf from } $1$ to $k$ {\bf do }\\
     \> \> \> $\displaystyle C_{N,k,m}:= \frac{N-1}{\lceil\frac{N-m}{k}\rceil}$\\
     \> \> \> $\displaystyle  V_{N,k,m}:= \sum_{m=1}^{\lceil\frac{N-m}{k}\rceil} |X(m+i\,k)-X(m+(i-1)\,k)|$. \\
     \> \> \> $\displaystyle  L_{m}(k):= \frac{1}{k} C_{N,k,m}V_{N,k,m}$. \\
     \> \> {\bf end} \\
     \>$L(k):=\frac{1}{k}\sum_{m=1}^{k}L_{m}(k)$\\
     \> {\bf end}\\
     \> $\mathcal{I}:=\{k \in \{1,..., k_{max}\} | L(k)\neq 0\}$\\
     \>  $\mathcal{Z}:=\{(\ln(\frac{1}{k}), \ln(L(k))) |k \in \mathcal{I}\}$\\
     \>{\bf if} $|\mathcal{I}|=1$ {\bf then}\\
     \>\> $D=1$\\
     \>{\bf else} \\
     \>\> $D$ is the slope of the best fitting (least square) \\
     \>\> affine function through $\mathcal{Z}$\\
     \>\> $HFD(X, N, k_{max})=D$\\
     \>{\bf end}\\
     \>  {\bf end}\\
\end{tabbing}}
\caption{Algorithm to compute the Higuchi fractal dimension for one dimensional series} \label{fig:algo_hig_1D}
\end{figure}

In the next examples we consider $X=[0,1]$ and the $\varepsilon$-net $\hat{X}=\{\frac{0}{M},\,\frac{1}{M},\, \ldots, \, \frac{M}{M}\}$, for a fixed number $M \geq 2$. All discretizations are made through a projection to the closer point in $\hat{X}$. So an idempotent measure $\mu= \bigoplus_{x \in X} \lambda(x) \odot \delta_x$ will be replaced by aa discrete version $\hat{\mu}= \bigoplus_{\hat{x} \in \hat{X}} \hat{\lambda}(\hat{x}) \odot \delta_{\hat{x}}$ which is identified to the series of values $\{\hat{\lambda}(\frac{0}{M}),\,\hat{\lambda}(\frac{1}{M}),\, \ldots, \, \hat{\lambda}(\frac{M}{M})\}$, those are the values showed in the drawings below.

\begin{example} \label{ex:example01 1-D draw}
    Let us consider $X=[0,1]$ endowed with the usual Euclidian induced metric. For each $j=1,2,\ldots$ we define the map $\phi_j: X \to X $ by
    $$\phi_j(x):= \frac{1+ x}{2^{j}}.$$
    We also consider the weights $q_j(x):= -(j-1)^2, \; j\geq 1.$

    Let $S=(X,\phi,q)$ be the associated mpCIFS, which is normalized because $\bigoplus_{j \in \mathbb{N}} -(j-1)^2 =-(1-1)^2=0$.

    We run the algorithm \texttt{DeterminIFSIdempMeasureDraw}($S$) with $M=1000$ and $N=30$ iterations.
\begin{figure}[H]
  \centering
  \includegraphics[width=5cm]{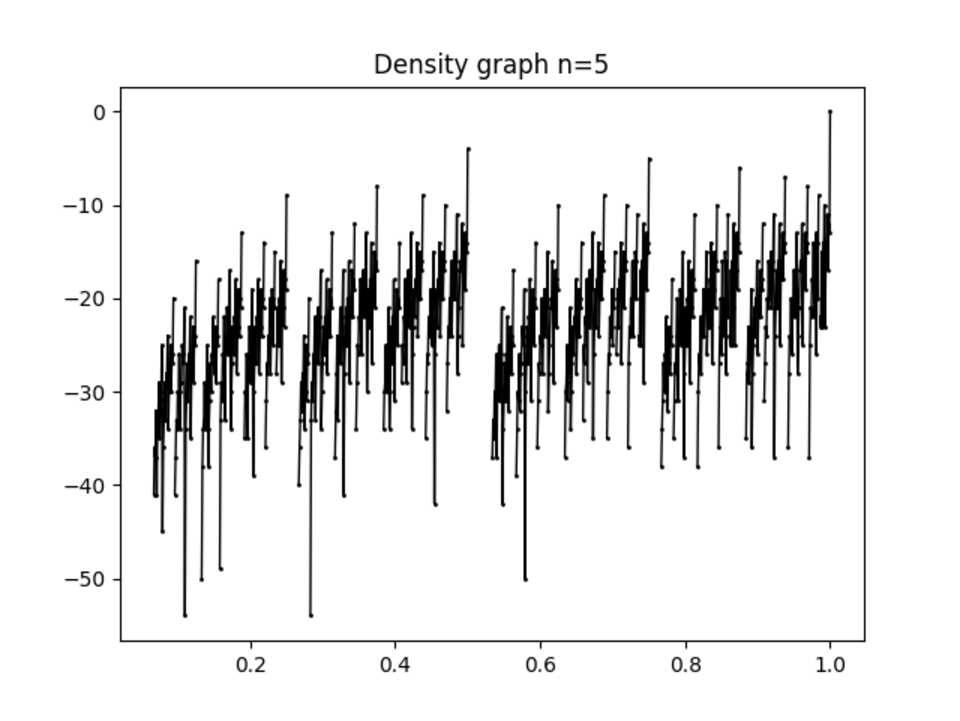} \includegraphics[width=5cm]{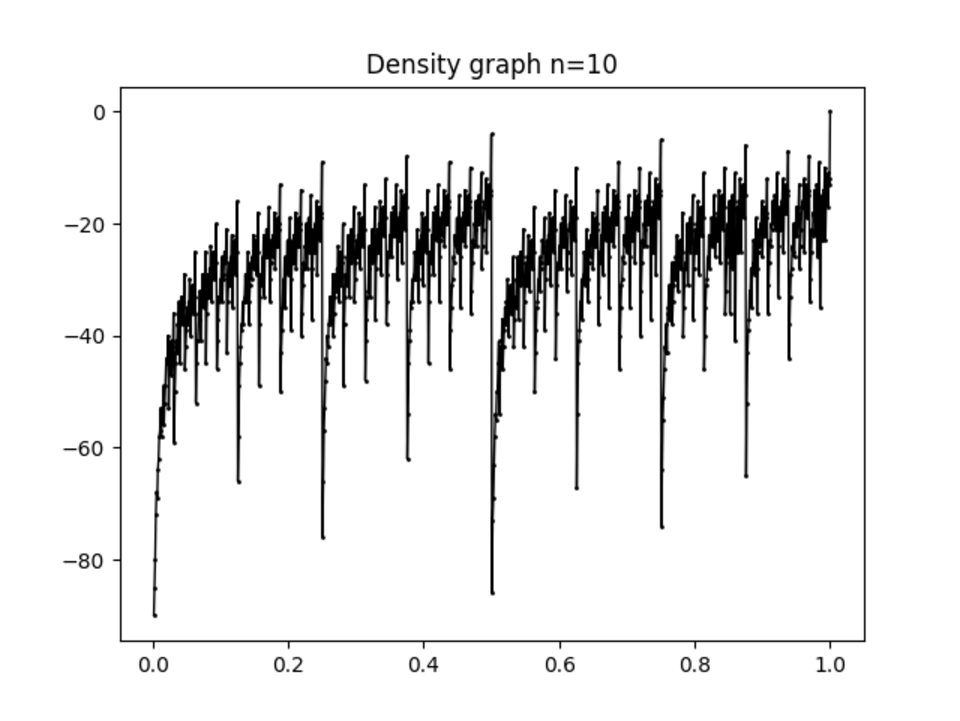}
  \includegraphics[width=5cm]{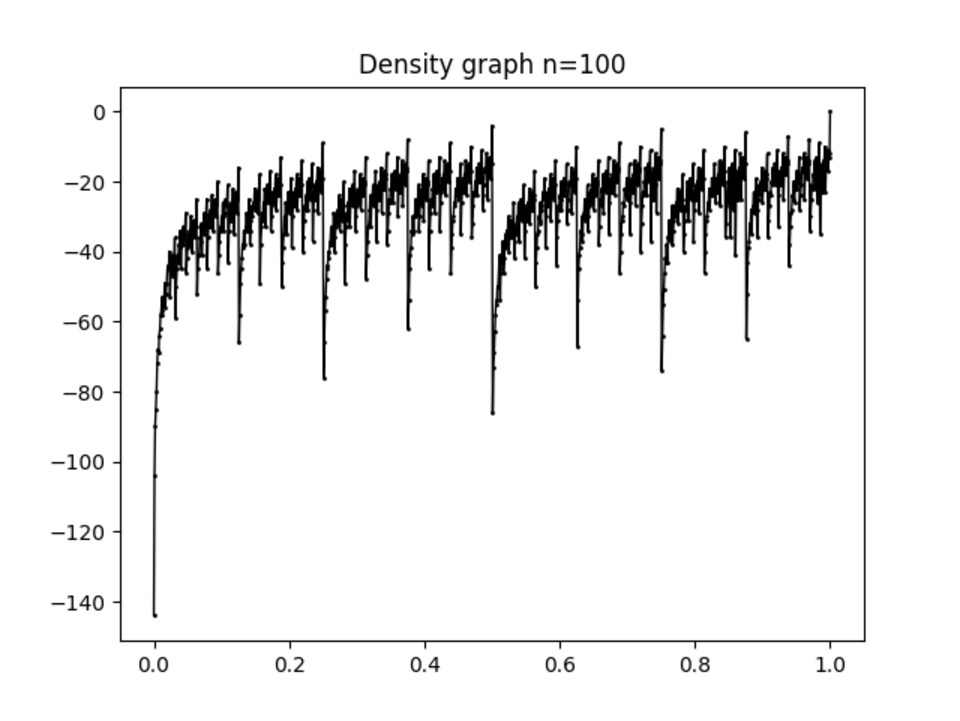}
  \caption{Approximation scheme for $\mu_{5}$, $\mu_{10}$ and $\mu_{100}$.}\label{fig:ex01 draw}
\end{figure}

Now we estimate the Higuchi fractal dimension of each approximation.
\begin{figure}[H]
  \centering
  \includegraphics[width=5cm]{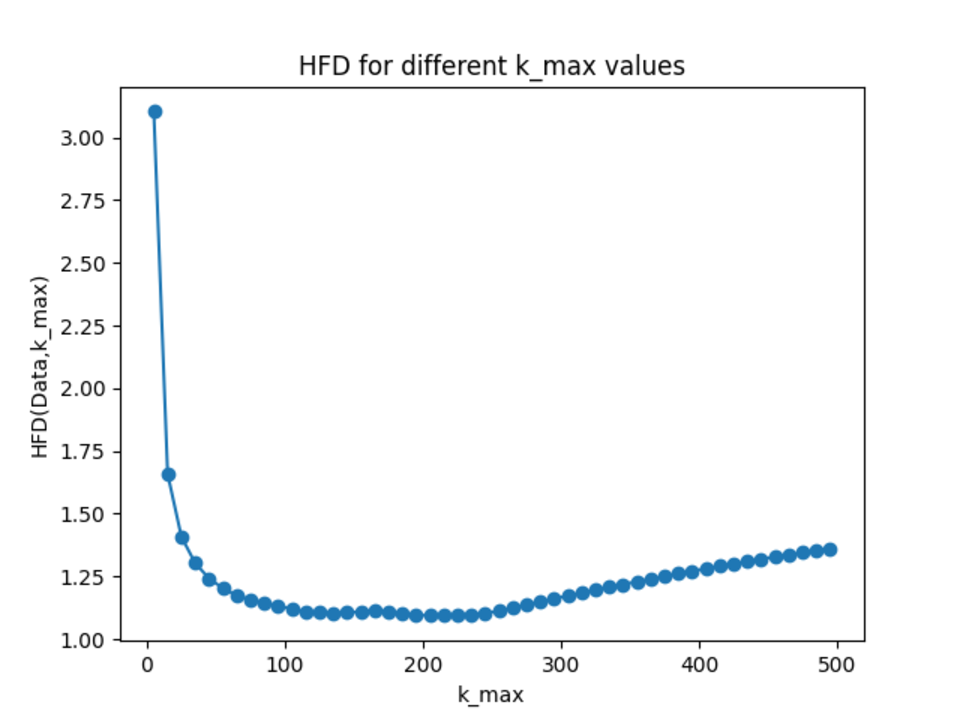} \includegraphics[width=5cm]{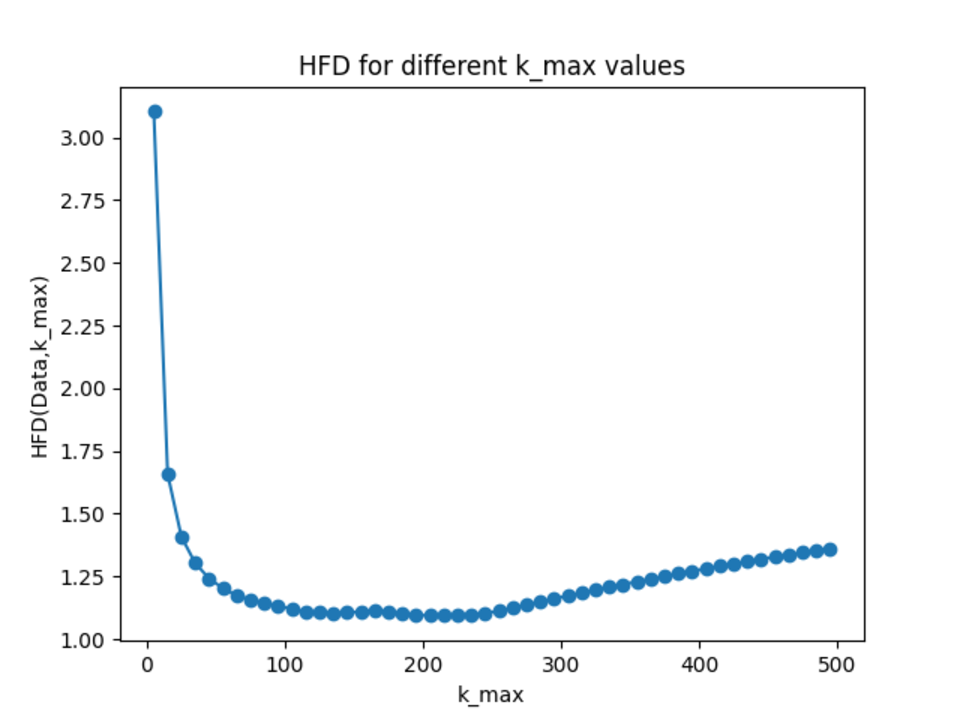}
  \includegraphics[width=5cm]{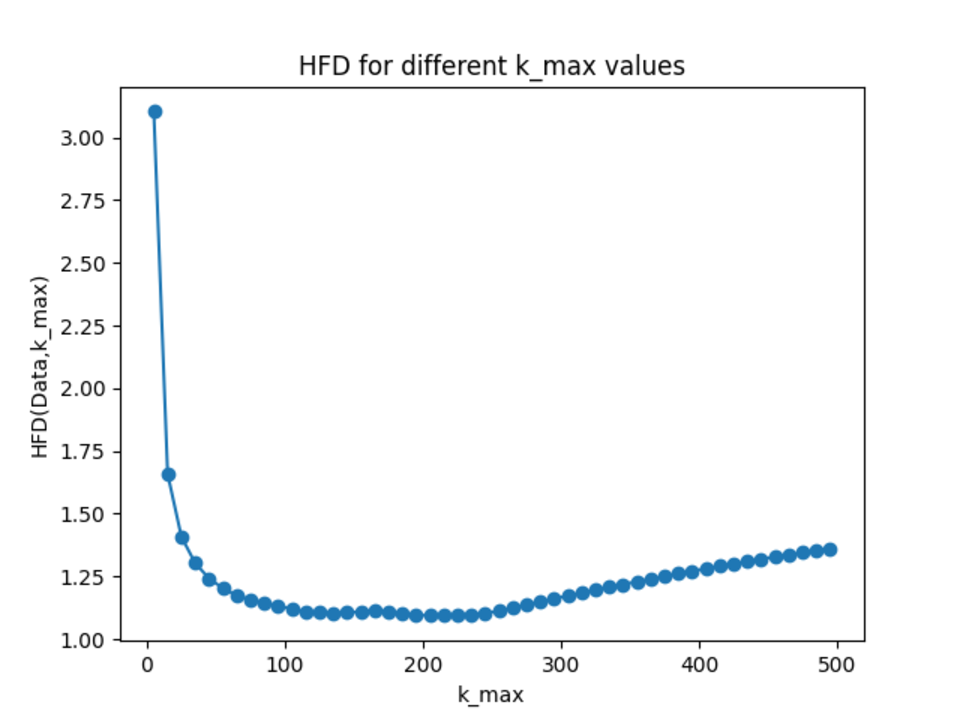}
  \caption{HFD for $\hat{\lambda}_{5}(x)$, $\hat{\lambda}_{10}(x)$ and $\hat{\lambda}_{100}(x)$.}\label{fig:ex01 HFD}
\end{figure}
As we can see, the value in which the fractal dimension estabilizes is around $k_{max}=200$ resulting in a fractal dimension of $1.09$.
\end{example}

 \begin{example}\label{ex:example02 1-D draw}
    Let us consider $X=[0,1]$ endowed with the usual Euclidian induced metric. For each $j=1,2,\ldots$ we define the map $\phi_j: X \to X $ by
    $$\phi_j(x):= \frac{1+ x}{2^{j}}.$$
    We also consider the weights
    $$q_j(x):= -\frac{1}{2^j}.$$
    Let $S=(X,\phi,q)$ be the associated mpCIFS, which is normalized because $\bigoplus_{j \in \mathbb{N}} -\frac{1}{2^j} =0$.

    We run the algorithm \texttt{DeterminIFSIdempMeasureDraw}($S$) with $M=1000$ and $N=30$ iterations.
\begin{figure}[H]
  \centering
  \includegraphics[width=5cm]{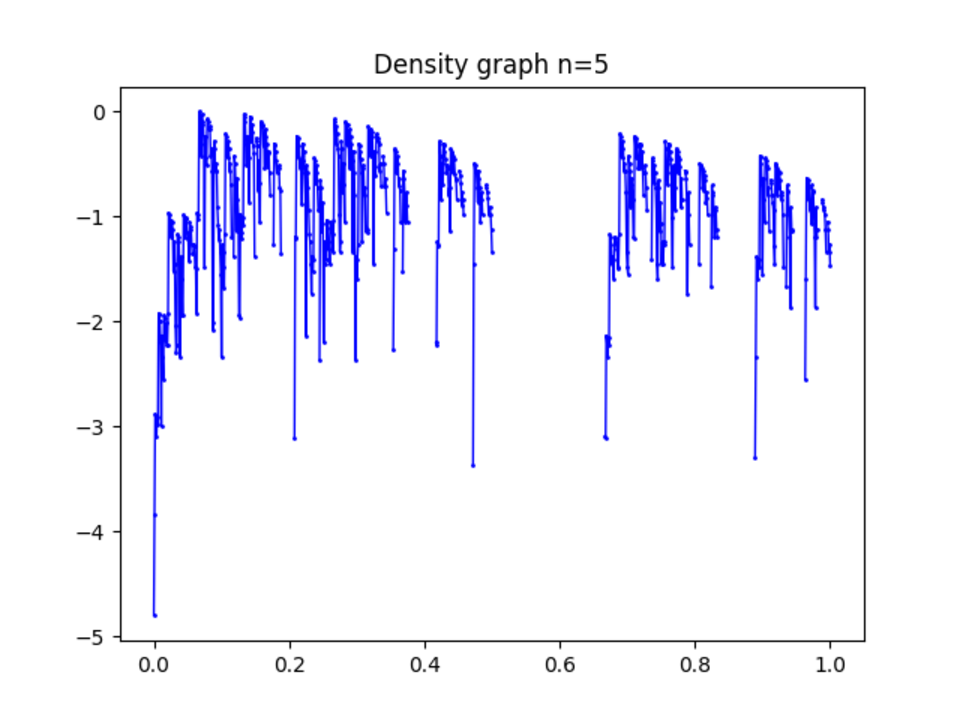} \includegraphics[width=5cm]{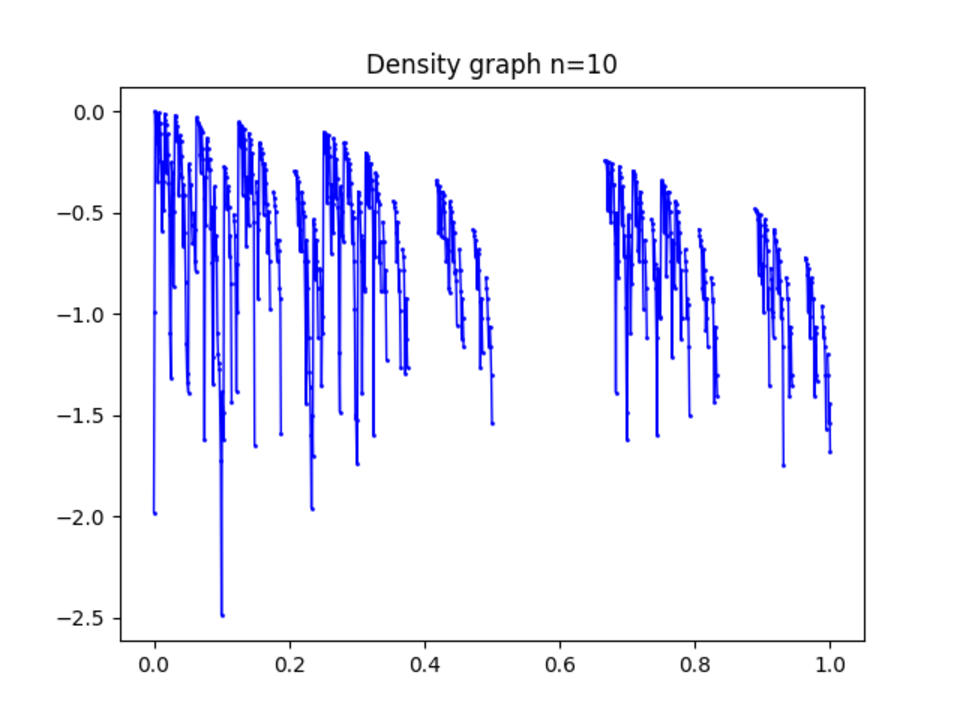}
  \includegraphics[width=5cm]{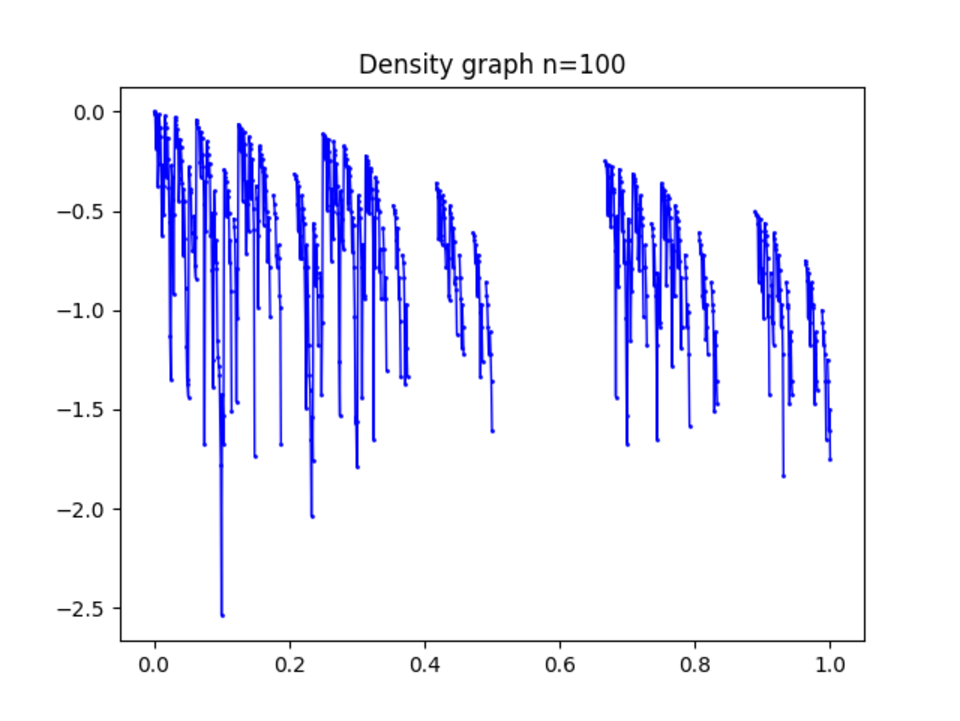}
  \caption{Approximation scheme for $\mu_{5}$, $\mu_{10}$ and $\mu_{100}$.}\label{fig:ex02 draw}
\end{figure}

Now we estimate the Higuchi fractal dimension of each approximation.
\begin{figure}[H]
  \centering
  \includegraphics[width=5cm]{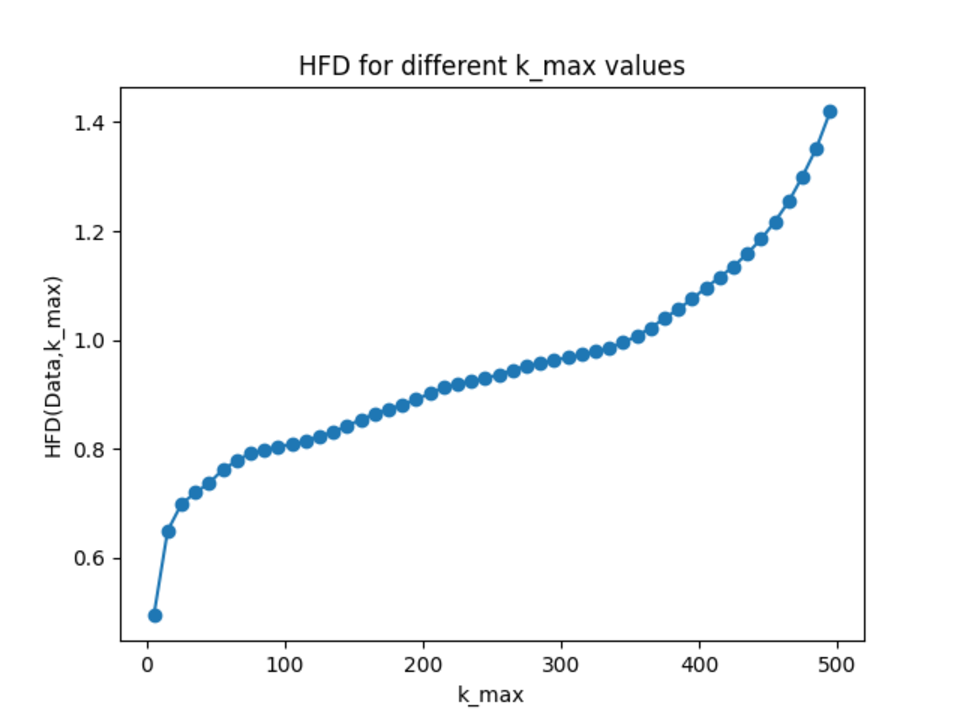} \includegraphics[width=5cm]{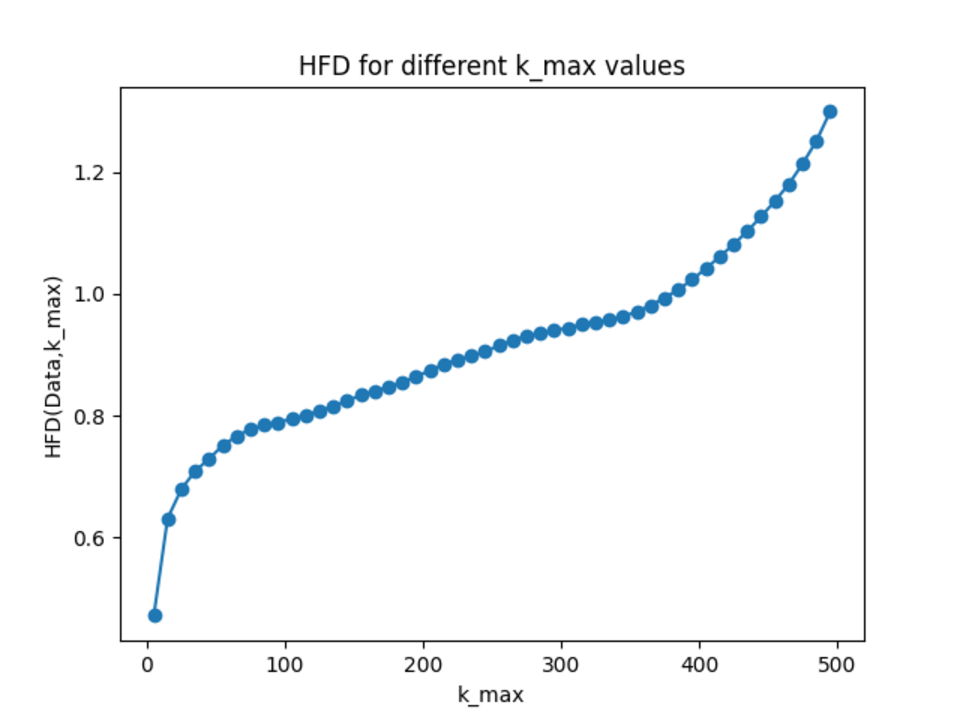}
  \includegraphics[width=5cm]{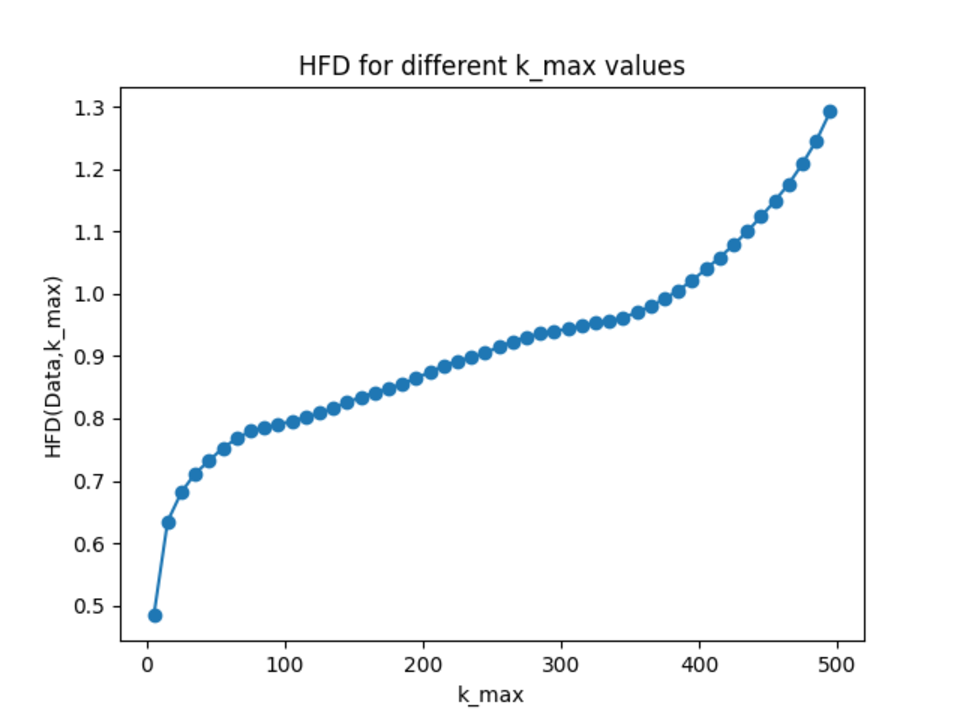}
  \caption{HFD for $\hat{\lambda}_{5}(x)$, $\hat{\lambda}_{10}(x)$ and $\hat{\lambda}_{100}(x)$.}\label{fig:ex02 HFD}
\end{figure}
As we can see, the value where the fractal dimension do not estabilizes attaining a fractal dimension of $1.29$ around $k_{max}=500$.
\end{example}

\subsection{2-D idempotent fractals and grey level images}\label{subsec:higuchi 2-d}

Our primary inspiration comes from Spasi\'{c}'s work \cite{Spa14}, which proposes a generalization of the Higuchi fractal dimension for bi-dimensional data by computing areas instead of lengths. The works of \cite{Aha11} and \cite{ASR15} are also noteworthy, as they introduce several 2D generalization algorithms, including the one in \cite{Spa14}, and evaluate the differences between these methods in computing the fractal dimension for digital images and other bi-dimensional data sets. However, it remains unclear whether there are significant advantages to using a specific model. The generalization idea in \cite{Spa14} relates to surface analysis in space and is based on Higuchi’s method for estimating fractal dimension values.

As highlighted by \cite{Sec01}, any compact subset of a metric space can be the attractor of a countable Iterated Function System (IFS). Therefore, our approach to 2-D figures is quite general regarding the inverse fractal problem. The choice of weights \( q_j \) will affect only the grey levels and thus the Higuchi dimension.

We now describe the approach in \cite{Spa14}. Given a time series \( X:\{1, \ldots, N\}^2 \to \mathbb{R} \), we aim to estimate the areas of \( X \) at different scales. Essentially, the data is organized as a matrix \( X = [x_{ij}] \), representing a discrete version of the function \( X \). We then consider triangles formed by \( x_{i,j}, x_{i+1,j}, \) and \( x_{i,j+1} \) on the graph of \( X \) at the coordinates \( (i,j), (i,j+1), (i+1,j) \). These triangles approximate the surface \( X \), and the area of \( X \) can be approximated by summing the areas of these triangles. Different dimensions are obtained by varying the scale of this procedure. The algorithm is detailed in \cite{Spa14} and \cite{ASR15}, and it adapts the procedure described in the original paper \cite{Hig88}:

\begin{figure}[H]
{\tt
\begin{tabbing}
aaa\=aaa\=aaa\=aaa\=aaa\=aaa\=aaa\=aaa\= \kill
    \>  \texttt{2D HFD algorithm} \\
    \> Input: Choose $2\leq N \in \mathbb{N}$, $X:\{1,...,N\}^2 \to \mathbb{R}$ and $2 \leq k_{max} \leq \lceil\frac{N}{2}\rceil$\\
     \> Output: Higuchi Fractal Dimension $HFD(X, N, k_{max})$\\
     \>  {\bf for } $k$ {\bf from } $1$ to $k_{max}$ {\bf do } \\
     \> \>   {\bf for } $n$ {\bf from } $1$ to $k$ {\bf do }\\
     \> \> \>  {\bf for } $m$ {\bf from } $1$ to $k$ {\bf do }\\
     \> \> \> \> $\displaystyle C_{N,k,n,m}:= \frac{N-1}{\lceil\frac{N-m}{k}\rceil} \, \frac{N-1}{\lceil\frac{N-n}{k}\rceil}$\\
     \> \> \> \> $\displaystyle  V_{N,k,m}:= \sum_{n=1}^{\lceil\frac{N-n}{k}\rceil}\sum_{m=1}^{\lceil\frac{N-m}{k}\rceil} |X(n+(i-1)\,k, m+j\,k)-X(n+(i-1)\,k, m+j\,k)|$ \\
     \> \> \> \> $|X(n+i\,k, m+j\,k)-X(n+(i-1)\,k, m+j\,k)|+$\\
     \> \> \> \> $|X(n+i\,k, m+j\,k)-X(n+i\,k, m+(j-1)\,k)|$ \\
     \> \> \> \> $|X(n+i\,k, m+(j-1)\,k)-X(n+(i-1)\,k, m+(j-1)\,k)|$. \\
     \> \> \> \> $\displaystyle  A_{n,m}(k):= \frac{1}{k^2} C_{N,k,n,m}V_{N,k,n,m}$. \\
     \> \> \> {\bf end} \\
     \> \> $A(k):=\frac{1}{2 k^2}\sum_{n=1}^{k}\sum_{m=1}^{k}A_{n,m}(k)$\\
     \> {\bf end}\\
     \> $\mathcal{I}:=\{k \in \{1,..., k_{max}\} | A(k)\neq 0\}$\\
     \> $\mathcal{Z}:=\{(\ln(\frac{1}{k^2}), \ln(A(k))) |k \in \mathcal{I}\}$\\
     \>{\bf if} $|\mathcal{I}|=1$ {\bf then}\\
     \>\> $D=1$ \\
     \>{\bf else} \\
     \>\> $D$ is the slope of the best fitting (least square) \\
     \>\> affine function through $\mathcal{Z}$\\
     \>\> $HFD(X, N, k_{max})=1+D$\\
     \>{\bf end}\\
     \>  {\bf end}\\
\end{tabbing}}
\caption{Algorithm to compute the Higuchi fractal dimension for two dimensional series} \label{fig:algo_hig_2D}
\end{figure}

\begin{example}\label{ex:example01 2-D draw}
    Let us consider $X=[0,1]^2$ endowed with the usual Euclidian induced metric. For each $j=1,2,\ldots$ we define the map $\phi_j: X \to X $ by
    $$\phi_j(x):=
    \left\{
      \begin{array}{ll}
        (\frac{x_1 }{2},\frac{x_2  }{2}), & j=4i \\
        (\frac{x_1 +1}{2},\frac{x_2  }{2}), & j=4i+1 \\
        (\frac{x_1 }{2},\frac{x_2 + 1}{2}), & j=4i+2 \\
        (\frac{x_1 + 1}{2},\frac{x_2 + 1}{2}), & j=4i+3. \\
      \end{array}
    \right.
$$
The fractal generated by the CIFS alone is a kind of ``checkered carpet".     We also consider the weights $q_j(x):= -(j-1)^2, \; j\geq 1.$

    Let $S=(X,\phi,q)$ be the associated mpCIFS, which is normalized because $\bigoplus_{j \in \mathbb{N}} -(j-1)^2 =-(1-1)^2=0$.

    We run the algorithm \texttt{DeterminIFSIdempMeasureDraw}($S$) with $M=256$ and $N=15$ iterations.
\begin{figure}[H]
  \centering
  \includegraphics[width=7cm]{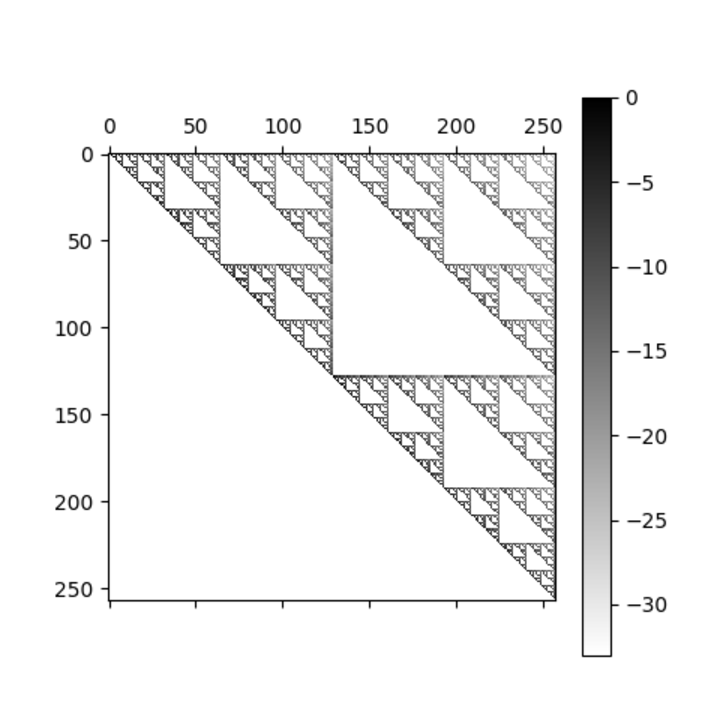} \includegraphics[width=7cm]{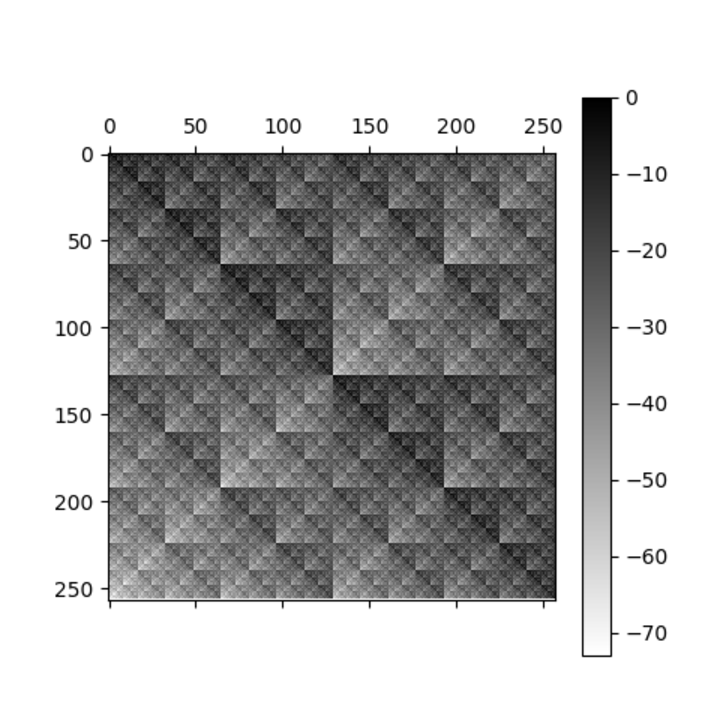}
  \caption{Approximation scheme for $\mu_{3}$ and $\mu_{15}$.}\label{fig:ex01 draw-2D}
\end{figure}

Now we estimate the Higuchi fractal dimension of each approximation.
\begin{figure}[H]
  \centering
  \includegraphics[width=7cm]{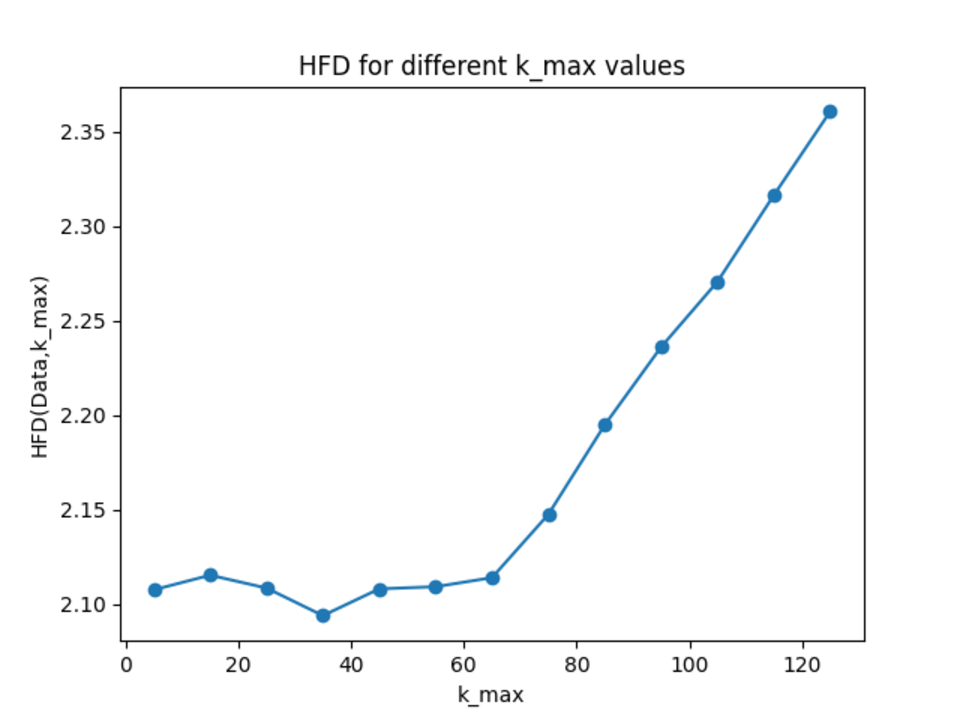} \includegraphics[width=7cm]{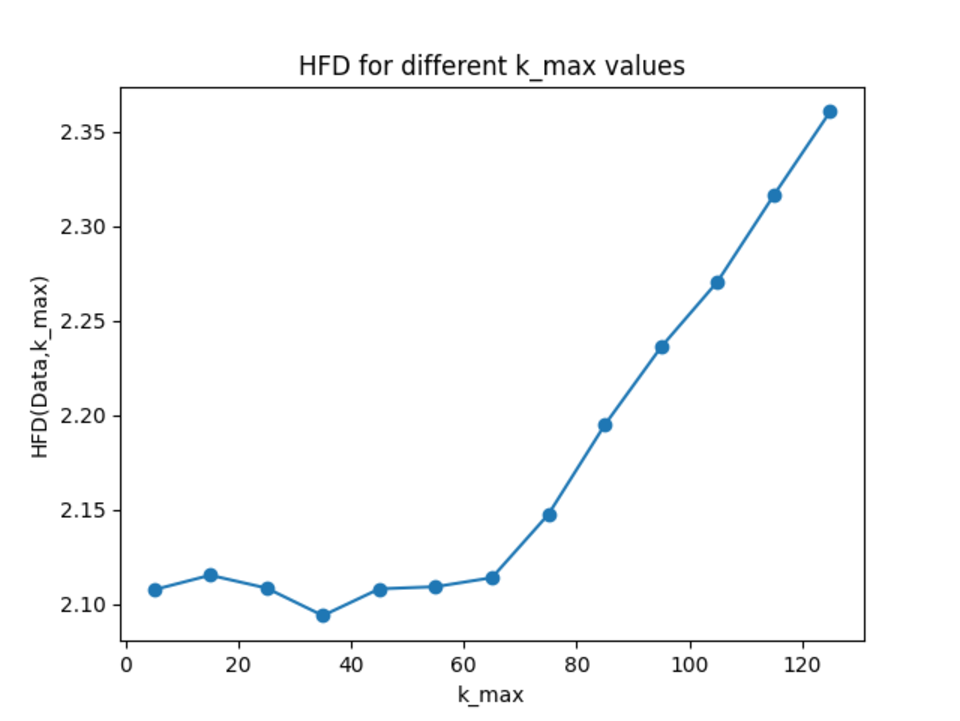}
  \caption{HFD for $\hat{\lambda}_{3}(x)$ and $\hat{\lambda}_{15}(x)$.}\label{fig:ex01 HFD-2D}
\end{figure}
As we can see, the value where the fractal dimension estabilizes is around $k_{max}=65$ resulting in a fractal dimension of $2.11$.
\end{example}

\begin{example}\label{ex:example02 2-D draw}
    Let us consider $X=[0,1]^2$ endowed with the usual Euclidian induced metric. For each $j=1,2,\ldots$ we define the map $\phi_j: X \to X $ by
    $$\phi_j(x):=
    \left\{
      \begin{array}{ll}
        (0.008\cdot x_{1}+0.1, 0.008\cdot x_{2}+0.04), & j=1 \\
        (0.5\cdot x_{1}+.25, 0.5\cdot x_{2}+.4), & j=2 \\
        (0.355\cdot x_{1}-.355\cdot x_{2}+.266, 0.355\cdot x_{1}+.355\cdot x_{2}+0.078), & j=3 \\
        (0.355\cdot x_{1}+0.355\cdot x_{2}+0.378\cdot (1-\frac{1}{j}),&\\ -0.355\cdot x_{1}+0.355\cdot x_{2}+0.434\cdot (1-\frac{1}{j})), & j \geq 4. \\
      \end{array}
    \right.
$$
The fractal generated by the CIFS alone is a kind of ``rotten" maple leaf.
    We also consider the weights
    $$q_j(x):= -\frac{1}{2^j}.$$
    Let $S=(X,\phi,q)$ be the associated mpCIFS, which is normalized because $\bigoplus_{j \in \mathbb{N}} -\frac{1}{2^j} =0$.

    We run the algorithm \texttt{DeterminIFSIdempMeasureDraw}($S$) with $M=256$ and $N=20$ iterations.
\begin{figure}[H]
  \centering
  \includegraphics[width=7cm]{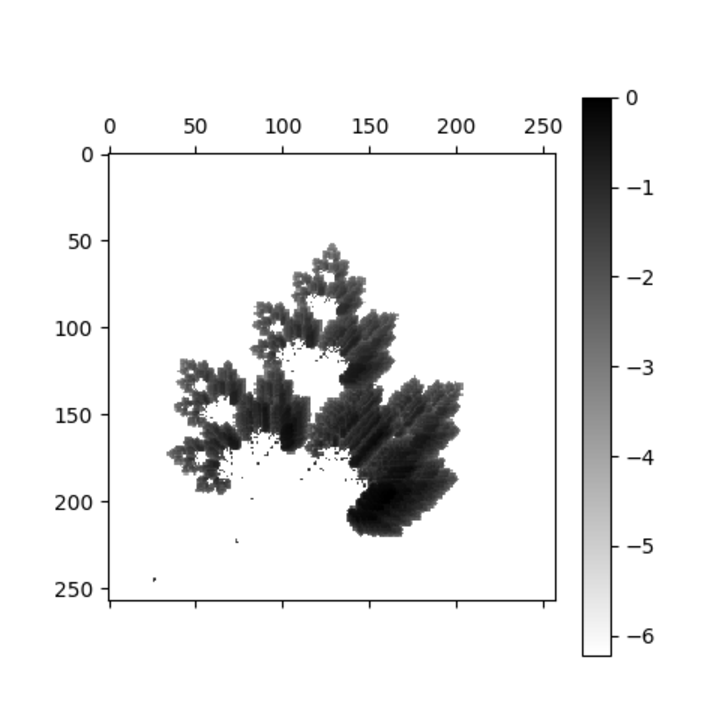} \includegraphics[width=7cm]{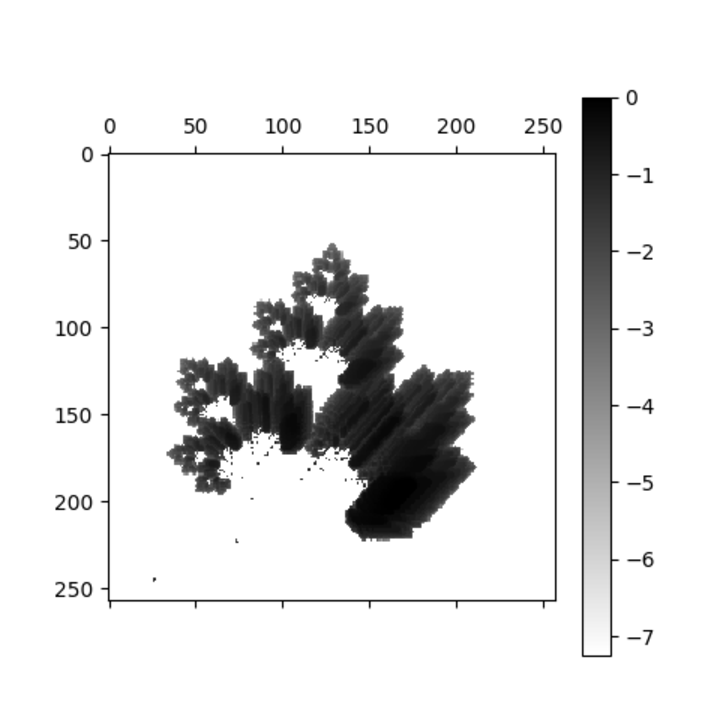}
  \caption{Approximation scheme for $\mu_{8}$ and $\mu_{15}$.}\label{fig:ex02 draw-2D}
\end{figure}

Now we estimate the Higuchi fractal dimension of each approximation.
\begin{figure}[H]
  \centering
  \includegraphics[width=7cm]{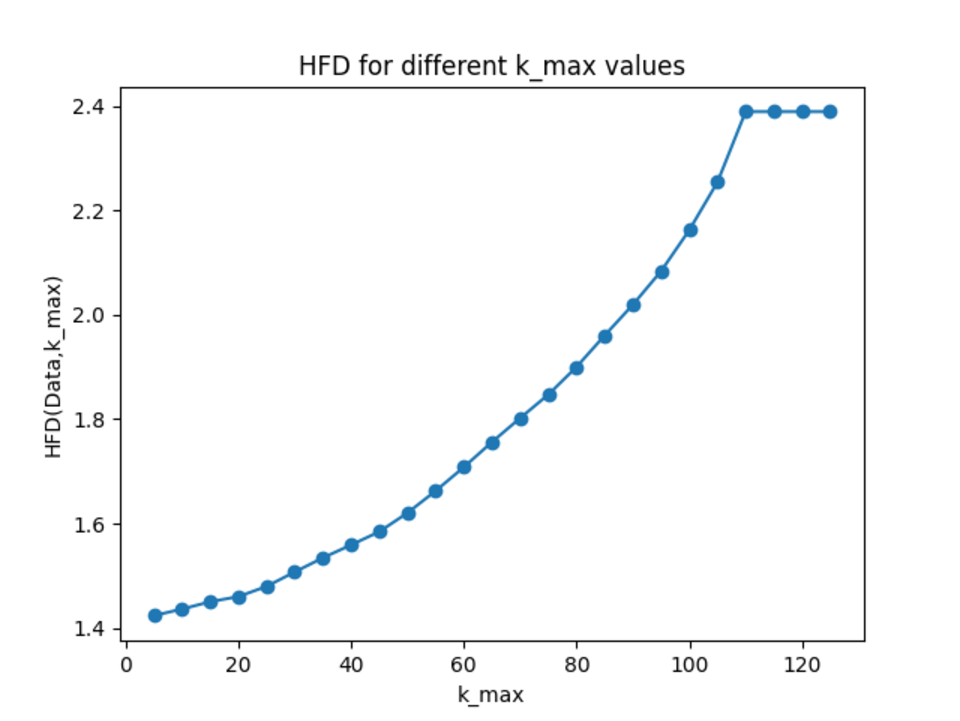} \includegraphics[width=7cm]{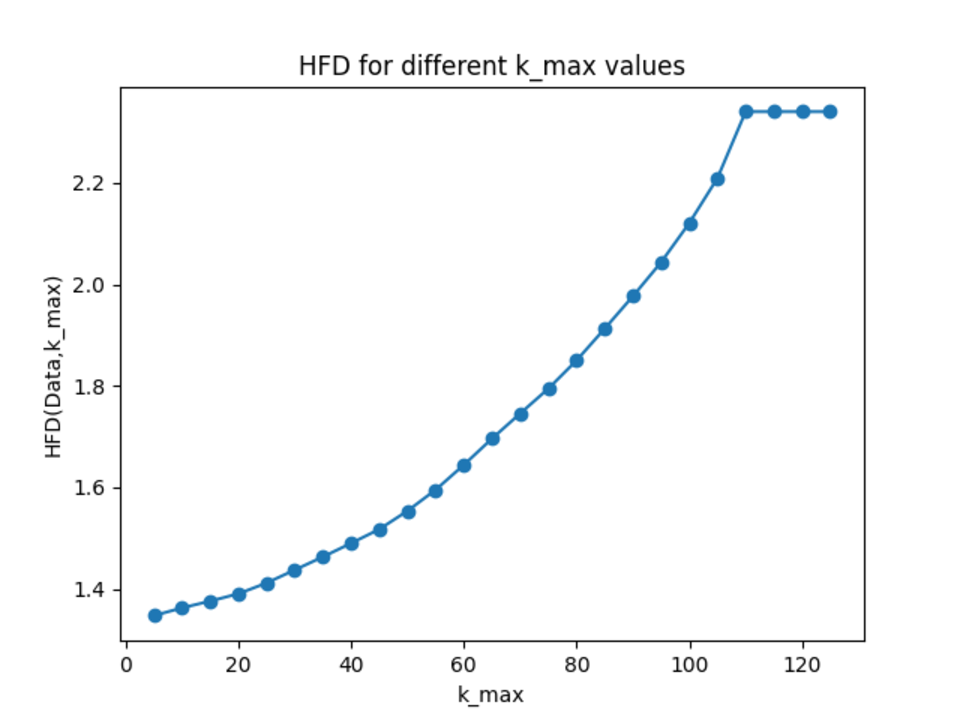}
  \caption{HFD for $\hat{\lambda}_{8}(x)$ and $\hat{\lambda}_{15}(x)$.}\label{fig:ex02 HFD-2D}
\end{figure}
As we can see, the value where the fractal dimension estabilizes is around $k_{max}=100$ resulting in a fractal dimension of $2.25$ for $j=8$ and  $2.34$ for $j=15$.
\end{example}

\textbf{Conclusion:} We defined the \((\delta, \varepsilon)\)-Higuchi Fractal Dimension applied to the fuzzyfication of discrete \(\delta\)-approximations, with respect to an \(\varepsilon\)-net, of the partial idempotent attractors of a max-plus contractive iterated function system (mpCIFS). Through various examples, we demonstrated that the Higuchi Fractal Dimension is an effective measure of the increasing complexity of the fractals as \(n\) approaches infinity, thereby approximating the actual attractor of the mpCIFS.

In these concrete examples, \(\delta\) represents the resolution number obtained when fixing the \(\varepsilon\)-net, the IFS, and the number of iterations, as described in \cite{dCOS21} and \cite{dCOS20fuzzy}.


\end{document}